\documentclass{svjour2}                    
\smartqed  
\usepackage{graphicx}
\usepackage[T1]{fontenc}
\usepackage{amsmath}
\usepackage{amssymb}

\newtheorem{examples}[example]{Examples}
\newtheorem{remarks}[remark]{Remarks}
\usepackage{color}

\def\a {\alpha}

\def\o {\omega}
\def\O {\Omega}
\def \veps {\varepsilon}
\def \vphi {\varphi}

\def\CL { \mathcal{L}}

\def\CP {\mathcal{P}}
\def\CN {\mathcal{N}}
\def\CI {\mathcal{I}}
\def\CF {\mathcal{F}}

\def\CL {\mathcal{L}}
\def\CB {\mathcal{B}}
\def\CS {\mathcal{S}}
\def\CR {\mathcal{R}}
\def\CT { \mathcal{T}}
\def\CA { \mathcal{A}}
\def\CE {\mathcal{E}}

\def\CK {\mathcal{K}}
\def\CR {\mathcal{R}}

\def\E {\mathbb{E}}

\def\R {\mathbb{R}}
\def\N {\mathbb{N}}
\def\I {\mathbb{I}}

\def\tte {\tilde{\theta}}
\def\ti {\tilde{i}}
\def\tnu {\tilde{\nu}}
\def\tbe {\tilde{\beta}}
\def\tlm {\tilde{\lambda}}
\def\tu {\tilde{u}}
\def\tLm {\tilde{\Lambda}}
\def\tb {\tilde{b}}

\def\cO {\widehat{\Omega}}
\def\cF {\widehat{\CF}}
\def\cP {\widehat{\CP}}
\def\cW {\widehat{W}}
\def \cCT {\widehat{\mathcal{T}}}

\def\cnun {\hat{\nu}^{(n)}}
\def\cten {\hat{\theta}^{(n)}}
\def\cin {\hat{i}^{(n)}}
\def\cben {\hat{\beta}^{(n)}}
\def\clmn {\hat{\lambda}^{(n)}}
\def\cLmn {\widehat{\Lambda}^{(n)}}
\def\cun {\hat{u}^{(n)}}

\def\cbe {\hat{\beta}}
\def\cnu {\hat{\nu}}
\def\cte {\hat{\theta}}
\def\ci {\hat{i}}

\def\ptm {\CP_{t,m}}

\begin{document}

\title{Probabilistic View of Explosion in an Inelastic Kac Model
}


\author{Andrea Bonomi         \and
        Eleonora Perversi \and
        Eugenio Regazzini 
}


\institute{A. Bonomi \at
              Dipartimento di Matematica, Universit\`a degli Studi di Pavia, 27100 Pavia, Italy \\
              \email{andrea.bonomi01@universitadipavia.it}           
           \and
           E. Perversi \at
              Dipartimento di Matematica, Universit\`a degli Studi di Pavia, 27100 Pavia, Italy \\
              \email{eleonora.perversi@unipv.it}           
           \and
           E. Regazzini \at
              Dipartimento di Matematica, Universit\`a degli Studi di Pavia, 27100 Pavia, Italy \\
              \email{eugenio.regazzini@unipv.it}\\
              Also affiliated to CNR-IMATI (Milano)
}

\date{Received: date / Accepted: date}

\maketitle

\begin{abstract}
Let $\{\mu(\cdot,t):t\geq0\}$ be the family of probability measures corresponding to the solution of the inelastic Kac model introduced in Pulvirenti and Toscani [\textit{J. Statist. Phys.} \textbf{114} (2004) 1453-1480]. It has been proved by Gabetta and Regazzini [\textit{J. Statist. Phys.} \textbf{147} (2012) 1007-1019] that the solution converges weakly to equilibrium if and only if a suitable symmetrized form of the initial data belongs to the standard domain of attraction of a specific stable law. In the present paper it is shown that, for initial data which are heavier-tailed than the aforementioned ones, the limiting distribution is improper in the sense that it has probability $1/2$ "adherent" to $-\infty$ and probability $1/2$ "adherent" to $+\infty$. It is explained in which sense this phenomenon is amenable to a sort of explosion, and the main result consists in an explicit expression of the rate of such an explosion. The presentation of these statements is preceded by a discussion about the necessity of the assumption under which their validity is proved. This gives the chance to make an adjustment to a portion of a proof contained in the above-mentioned paper by Gabetta and Regazzini.
\end{abstract}

\section{Introduction}\label{sec:intro}
This paper deals with the \textit{explosion} of the solution of the \textit{inelastic Kac model} introduced in \cite{PulvirentiToscani}, also studied, e.g., in \cite{BaLaRe,GabettaRegazzini2012} and extended to various fields in \cite{BaLa,BassLadMatth10,BaLaTo,BassettiPerversi,BassettiToscani,MatthesToscani,MaTo}. This model describes the evolution in time $t$ of the probability distribution (p.d., for short), say $\mu(\cdot,t)$, of the velocity of a particle in a granular gas subject to dissipative collisions and only for space independent data. More specifically, when two particles, with velocities $v$, $w$ respectively, collide, they change their velocities to $v^*$, $w^*$ given by 
\[
v^*=v c_p(\theta)-w s_p(\theta),\qquad w^*=w s_p(\theta)+v c_p(\theta)
\]
where $\theta$ is any angle in $(0,2\pi]$, $p$ is a nonnegative parameter, $c_p$ and $s_p$ are defined by
\[
c_p(\theta):=\cos\theta|\cos\theta|^p, \quad s_p(\theta):=\sin\theta|\sin\theta|^p.
\]
Written in terms of the Fourier-Stieltjes transform (\textit{characteristic function}, c.f. for short) $\vphi(\xi,t):=\int_\R e^{i\xi x}\mu(dx,t)$ for every $\xi$ in $\R$, the aforesaid model reduces to the equation

\begin{equation}\label{eqCauchy}
	\left\{	\begin{aligned}
					& \frac{\partial}{\partial t}\varphi(\xi,t)=\frac{1}{2\pi}\int_{0}^{2\pi}\varphi(\xi c_p(\theta),t)\varphi(\xi s_p(\theta),t)d\theta-\varphi(\xi,t)\qquad(t>0) \\
					& \vphi(\xi,0^+)=\vphi_0(\xi)\\
					\end{aligned}
	\right.
\end{equation}
where $\vphi_0$ stands for the c.f. associated with the initial velocity p.d. $\mu_0$. The parameter $p$ can be viewed as an \textit{index of inelasticity} in the model, perfect elasticity being realized for $p=0$, i.e. when \eqref{eqCauchy} reduces to the more renowned Kac's equation. Motivations for the study of dissipative systems ($p>0$) are given in \cite{PulvirentiToscani} and, for more realistic multidimensional equations, in many papers among which  \cite{BolleyCarrillo,BobylevCarrilloGamba,BobylevCercignani,CarrilloCordierToscani}. Interesting reviews can be found in \cite{CarrilloToscani2007,Villani2002,Villani2006}. Coming back to \eqref{eqCauchy}, it is well-known that it has a solution, which is unique within the class of the c.f.'s associated with all probability measures (p.m.'s, for short) on $\R$. 

There are classical problems in the kinetic theory of gases which lead to the study of the long-time behaviour of the solutions of the equations of interest. Apropos of \eqref{eqCauchy}, answers to this kind of problems have been efficaciously expressed by means of the \textit{probability distribution function} (p.d.f., for short) of the initial datum $\mu_0$, i.e. $F_0(x):=\mu_0((-\infty,x])$ for every $x$ in $\R$. More specifically, by resorting to the symmetric p.d.f. defined by
\[
F^*_0(x):=\frac{1}{2}\Big[F_0(x)+1-F_0(-x)\Big]
\]
at every continuity point $x$, Theorem 1 in \cite{GabettaRegazzini2012} states that the validity of the condition

\begin{equation}\label{NDA}
 \lim_{x\to+\infty}x^\a[1-F_0^*(x)]=c_0
\end{equation}
$-$ where $c_0$ is a nonnegative real number and $\a$ depends on $p$ according to $\a=2/(1+p)$ $-$ implies that $\mu(\cdot,t)$ converges weakly, as $t\to+\infty$, to a probability measure $\mu_\infty$ with c.f.
\[
\vphi_\infty(\xi)=e^{-a_0|\xi|^\a}\qquad(\xi\in\R)
\]
where 
\[
 a_0=2c_0 \lim_{T\to+\infty} \int_0^T \frac{\sin x}{x^\a}dx.
\]
This shows that the limit of $\mu(\cdot,t)$ is, for $c_0>0$, a symmetric stable law, which reduces to the unit mass $\delta_0$ if $c_0= 0$, i.e., if the tail $(1-F_0^*(x))$ is of smaller order (as $x\to+\infty$) than the heavy tail of the initial data attracted to a nondegenerate $\mu_\infty$ ($c_0>0$). Of course, this happens if $\int_\R|x|^d dF^*_0(x)<+\infty$ for some $d\geq \a$. Moreover, Theorem 1 in \cite{GabettaRegazzini2012} states that \eqref{NDA} is also necessary for  convergence of the solution. This way, we see that $\mu(\cdot,t)$ is not weakly convergent if and only if either $\liminf_{x\to+\infty}x^\a\Big[1-F_0^*(x)\Big]<\limsup_{x\to+\infty}x^\a\Big[1-F_0^*(x)\Big]$ or $\lim_{x\to+\infty}x^\a\Big[1-F_0^*(x)\Big]=+\infty$. In the present paper we aim at investigating into the long-time behaviour of the solution of \eqref{eqCauchy} when the initial datum meets $\limsup_{x\to+\infty}x^\a\Big[1-F_0^*(x)\Big]=+\infty$, i.e., the prior p.d. is, in a sense, \textit{ultra-heavy-tailed}. Apropos of this terminology, we have just noted that the reaching of nondegenerate equilibria in a dissipative model is assured only by the adoption of heavy-tailed $\mu_0$ when $p>0$. Moreover, as a consequence of next Proposition \ref{prop2}, we can state that the choice of symmetrized initial data exhibiting tails of greater order than $1/x^\a$ $(x\to+\infty)$ is necessary to observe the explosion of the solution. These a few remarks enable us to explain the meaning of the term explosion. From a purely mathematical standpoint, it is here used to mean that equalities $\lim_{t\to+\infty}\mu((-\infty,-a],t)=\lim_{t\to+\infty}\mu([a,+\infty),t)=1/2$ hold \textit{for every} $a>0$. A physical interpretation can be given following the same argument as in Section 1 of \cite{CarGabReg}, where the phenomenon is analysed for Kac's equation ($p=0$). One notes therein that the equations of ultimate interest ought to be spatially inhomogeneous versions of \eqref{eqCauchy}. They describe the evolution of the joint distribution of the pair (velocity, position) of a particle. In such a spatially dependent setting, the homogeneous equation studied here can be viewed as a picture of the evolution of a conditional p.d. for the velocity, given a specific position $x$, in the presence of an initial conditional p.d. $\nu_{0,x}=\mu_0$ deduced from the initial joint law $\nu_0$ of the above pair on the phase space $\R^2$. What is important is to note that joint initial data $\nu_0$, producing proper (tight) limiting distributions, can be consistent with conditional p.d.'s $\nu_{0,x}$ that, at certain positions $x$ called "hot positions", exhibit ultra-heavy-tails. In this perspective, our next Theorem \ref{teorema} states that all the molecules involved in a "hot position" $x$ quickly pick up very high velocity, isotropically distributed, which then explode away from $x$. Whence, the use of the term explosion here is literal and physical, and it is not to be confused with its common use in the theory of stochastic processes. These remarks justify, in our opinion, the study of the spatially homogeneous equation for initial data having ultra-heavy-tails. Moreover, they give reasons for the finding of a \textit{rate of explosion} meant as maximal speed of divergence to infinity, as $t\to+\infty$, of those $a_t$'s for which $\lim_{t\to+\infty}\mu((-\infty,- a_t],t)=\lim_{t\to+\infty}\mu([a_t,+\infty),t)=1/2$. A lower bound for such a kind of rate of explosion will be given in Section \ref{sec:main}. 

As to the organization of the paper, the novel results above briefly described will be stated in Section \ref{sec:main} in an autonomous way, in the sense that they can be understood on the sole basis of the present introduction. In spite of this, we have decided to interpose, between Section \ref{sec:intro} and Section \ref{sec:main}, the presentation of some preliminary statements, which are essential for the proof,  deferred to Section \ref{sec:prove}, of our main theorem. Proofs of these preliminary statements and of other three preparatory propositions, formulated in Sections \ref{sec:main} and \ref{sec:prove}, are presented in Appendices A, B, C and E at the end of the paper. Appendix D contains some technical details to determine examples of upper bounds for $a_t$, given at the end of Section \ref{sec:main}. We take advantage of the wording of Proposition \ref{prop1} to make a correction to the argument used in the first part of Step 2 in Section \ref{sec:main} of \cite{GabettaRegazzini2012}.

\section{Preliminaries}\label{sec:preliminaries}
The proof of the main theorem to be formulated in the next section rests on a probabilistic scheme originally given in \cite{McKean,McKean67} and exploited, for example, in \cite{BaLa,BassLadMatth10,CarlenCarvalhoGabetta,CarlenCarvalhoGabetta2005,DoleraGabettaRegazzini,DoleraRegazzini2010,DoleraRegazzini2012,GabReg2006,GabettaRegazzini2008,GabettaRegazzini2010}. Hence, we first touch on the basics of such a scheme. A probability space $(\O,\CF,\CP_t)$ is defined for each $t>0$. Random elements $X$, $\tte$, $\ti$, $\tnu$ are defined on $(\O,\CF)$ so that they turn out to be \textit{stochastically independent} with respect to each $\CP_t$. Definitions and additional distributional properties, with respect to each $\CP_t$, are: 
\begin{itemize}
\item $X=(X_n)_{n\geq1}$ is a sequence of independent and identically distributed (i.i.d., for short) random numbers with common p.d. $\mu_0$.
\item $\tte=(\tte_n)_{n\geq1}$ is a sequence of i.i.d. random numbers uniformly distributed on $(0,2\pi]$.
\item $\ti=(\ti_n)_{n\geq1}$ is a sequence of independent integer-valued random numbers, with $\ti_n$ uniformly distributed on $\{1,\dots,n\}$.
\item $\tnu$ is an integer-valued random number such that $\CP_t\{\tnu=n\}=e^{-t}(1-e^{-t})^{n-1}$ for $n=1,2,\dots\;.$
\end{itemize}
In this probabilistic setting, one can state (see, e.g., Theorem 3 in \cite{GabettaRegazzini2008}) that the solution $\mu(\cdot,t)$ of \eqref{eqCauchy} is the p.d., under $\CP_t$, of the random number

\begin{equation}\label{V_t}
 V:=\sum_{j\geq1}\tbe_{j,\nu}X_j=(\tbe_{\nu},X)
\end{equation}
where $\tbe_{k}$ is defined recursively, for $k=1,2,\dots$, as follows. Firstly, one puts $\tbe_1:=(1,0,0,\dots)$ and then, for any $k\geq1$, 
\begin{equation*}
\begin{split}
\tbe_{k+1}&=(\tbe_{1,k+1},\dots,\tbe_{k+1,k+1},0,0,\dots)\\
&:=(\tbe_{1,k},\dots,\tbe_{\ti_k-1,k},c_p(\tte_k)\tbe_{\ti_k,k},s_p(\tte_k)\tbe_{\ti_k,k},\tbe_{\ti_k+1,k},\dots,\tbe_{k,k},0,0,\dots).
\end{split}
\end{equation*}
It should be noted that $(\ref{V_t})$ entails

\begin{equation*}
\begin{split}
\vphi(\xi,t)=\int_{\Omega}\prod_{j=1}^{\tnu(\omega)}\vphi_0(\tbe_{j,\tnu(\omega)}(\omega)\xi)\CP_t(d\omega)
\end{split}
\end{equation*}
which provides a \textit{fibering of the solution of \eqref{eqCauchy} into components which are c.f.'s of weighted sums of i.i.d. random numbers with common p.d. $\mu_0$}.  Thanks to this kind of fibering the study of problem \eqref{eqCauchy} can be associated with the central limit problem of the probability theory. It is useful to emphasise another representation, according to which 

\begin{equation*}
\begin{split}
\vphi(\xi,t)&=\sum_{n\geq1}e^{-t}(1-e^{-t})^{n-1}\int_{\Omega}\prod_{j=1}^{n}\vphi_0(\tbe_{j,n}(\omega)\xi)\CP_t(d\omega)\\
&= i e^{-t} \Big(\Im \vphi_0(\xi)\Big)+ \sum_{n\geq1}e^{-t}(1-e^{-t})^{n-1}\int_{\Omega}\prod_{j=1}^{n}\Re\Big(\vphi_0(\tbe_{j,n}(\omega)\xi)\Big)\CP_t(d\omega)
\end{split}
\end{equation*}
where $\Re z$ ($\Im z$, respectively) denotes the real (imaginary, respectively) part of a complex number $z$. This representation explains that, with a view to the study of the asymptotic behaviour of the solution of \eqref{eqCauchy}, one can confine oneself to considering initial data equal to $\Re \vphi_0$, that is, equal to the c.f. $\vphi^*_0$ associated with $F^*_0$. See the end of Subsection 1.1 in \cite{BaLaRe}.

In some passages of the paper we shall resort to the \textit{Skorokhod representation}. Then, we conclude this preliminary section by describing the use of such a representation to deduce useful properties of the coefficients $\tbe_{j,\tnu}$. Following the argument originally explained in \cite{ForLadReg}  and more recently applied in \cite{BaLaRe,DoleraRegazzini2012,GabettaRegazzini2012}, we introduce a random vector $W$ which displays the quantities that are used to characterize convergence in general forms of the central limit theorem. More specifically, for each $\o$ in $\O$ we set 
\[
W=W(\omega):=(\tnu(\o),\tte(\o),\ti(\o),\tbe(\o),\tlm(\o),\tLm(\o),\tu(\o))
\]
where $\tnu$, $\tte$, $\ti$ are the same as at the beginning of this section, and
\begin{itemize}
\item $\tbe=\tbe(\o)$ is the matrix whose $k$-th row is $\tbe_k$, $k\geq1$.
\item $\tlm=\tlm(\o):=(\tlm_1(\o),\dots,\tlm_{\tnu(\o)}(\o),\delta_0,\delta_0,\dots)$ and, for each $j$ in $\{1,\dots,$ $\tnu(\o)\}$, $\tlm_j(\o)$ is the p.d. determined by the c.f. $\xi\mapsto\vphi^*_0(|\tbe_{j,\tnu(\o)}(\o)|\xi)$, $\xi\in\R$.
\item $\tLm=$convolution of the elements of $\tlm$.
\item $\tu:=(\tu_k)_{k\geq1}$, with $\tu_k=\max_{1\leq j\leq\tnu}\tlm_j([-\frac{1}{k},\frac{1}{k}]^c)$ for every $k\geq1$.
\end{itemize}
The range of $W$ is a subset of
\[
S:=\overline{\N}\times [0,2\pi]^\infty\times \overline{\N}^\infty\times [0,1]^\infty\times (\CP(\overline{\R}))^\infty\times \CP(\overline{\R})\times [0,1]^\infty
\]
where, for any metric space $M$, $\CP(M)$ has to be understood as the set of all p.d.'s on the Borel class $\CB(M)$. $\CP(\overline{\R})$ is metrized consistently with the topology of weak convergence of p.m.'s and, thus, it turns out to be a separable, compact and complete metric space. Thus, $S$ is seen as separable, compact and complete metric space (Theorems 6.2, 6.4 and 6.5 in Chapter 2 of \cite{Parthasarathy}). Then, the family of p.m.'s $\{\CP_tW^{-1}:t\geq0\}$ is \textit{uniformly tight} on $\CB(S)$, which implies that any subsequence of this family contains a weakly convergent subsequence $Q_n:=\CP_{t_{n}}W^{-1}$ with $0\leq t_1<t_2<\dots$ and $t_n \nearrow +\infty$. Thus, we are in a position to apply Skorokhod's representation theorem (see, e.g., Theorem 6.7 in \cite{Billingsley}) according to which there are a probability space $\Big(\cO,\cF,\cP\Big)$ and a sequence of $S$-valued random elements 
\[
\cW_n:=(\cnun,\cten,\cin,\cben,\clmn,\cLmn,\cun),\qquad n\geq1
\]
defined on $\cO$ such that:
\begin{itemize}
\item For every $n$, the p.d. of $\cW_n$ is $Q_n$.
\item $\cW_n$ converges pointwise to a random element $\cW$ whose p.d. is the weak limit of $(Q_n)_{n\geq1}$.
\end{itemize}
In view of the former of these two properties, the equalities
\begin{equation}\label{ricorrenzabeta}
\begin{array}{ll}
&\cben_1=(1,0,0,\dots)\\
& \cben_{k+1}=(\cben_{1,k},\dots,\cben_{\cin_k-1,k},c_p(\cten_k)\tbe_{\cin_k,k},s_p(\cten_k)\cben_{\cin_k,k},\cben_{\cin_k+1,k},\\
&\qquad\qquad\dots,\cben_{k,k},0,0,\dots)
\end{array}
\end{equation}
hold for every $k$ and $n$, with $\cP$-probability $1$. This paves the way for a statement that will be used later in conjunction with an application of the central limit theorem to the convolution $\cLmn=\clmn_1\ast\clmn_2\ast\dots$.

\begin{proposition}\label{beta}
Let $(x_m)_{m\geq1}$ be a strictly increasing and divergent sequence of positive real numbers. Then, for every strictly positive $\veps$, there are a strictly increasing and divergent $\N$-valued sequence $(m_n)_{n\geq1}$ and a point $\o^*$ in $\cO$ $-$ which depend on $(x_m)_{m\geq1}$ and $\veps$ $-$ such that
\begin{equation}\label{limB}
\frac{1}{x_{m_K}}\leq \Big| \cben_{j,\cnu^{(K)}(\o^*)}(\o^*)\Big| \leq \frac{1}{x_{m_K}-\veps}\quad (n\geq1,\; K=1,\dots,n)
\end{equation}
holds for every $j=2i-1$ $\Big(i=1,\dots,\left\lfloor \frac{\cnu^{(K)}(\o^*)+1}{2}\right\rfloor\Big)$. Moreover, there is $a>0$, independent of $n$, for which
\begin{equation}\label{sumB}
\sum_{i=1}^{\lfloor \frac{\cnu^{(K)}(\o^*)+1}{2}\rfloor}| \cben_{2i-1,\cnu^{(K)}(\o^*)}(\o^*)|^\a\geq a \qquad (n\geq1,\; K=1,\dots,n)
\end{equation}
$\a$ being the same as in \eqref{NDA}.
\end{proposition}
The proof is deferred to Appendix A. In preparation of such a proof, it is worth hinting at an equivalent construction of the $\tbe$'s based on the so-called \textit{McKean trees}. One starts from the "root node" assuming that it has two "children": a "left child" and a "right child". With reference to the previous notation, this first step is tantamount to saying that $\ti_1=1$ with probability $1$. This step made, the two children are labelled $1,2$ following a left-to-right order. One of them, say $\ti_2$, is chosen "at random" and it is replaced by a couple of "children" giving raise to a tree with $3$ leaves, labelled by $1,2,3$ from left-to-right. One of them ($\ti_3$ in the usual notation) is selected "at random" and substituted, like above, with a couple of "children" to obtain a tree with $4$ leaves, and so on. A sequence of $(n-1)$ steps of this kind produces a specific McKean tree with $n$ leaves. If, at each step, an angle $\tte_j$ is chosen at random and $c_p(\tte_j)$ ($s_p(\tte_j)$, respectively) is associated with the "left child" ("right child", respectively) of $\ti_j$, for every $j$, then $\tbe_{k,n}$ turns out to be the product of all the $c_p(\cdot)$'s and $s_p(\cdot)$'s along the path that joins the $k$-th leaf and the "root node" in the tree. The number of steps (i.e. the number of factors of $\tbe_{k,n}$) in such a path is called \textit{depth} of the $k$-th leaf. Trees with $n=2^m$ leaves are called \textit{complete} if the depth of each leaf is constant ($=m$).

\section{Results}\label{sec:main}
As anticipated in Section 1, we here present a quantification of the rate of explosion of the solution of \eqref{eqCauchy} when $p>0$, i.e. $0<\a<2$. The case in which $p=0$ ($\a=2$) has been dealt with in \cite{CarGabReg}. In order that explosion occur, it is necessary that the solution of \eqref{eqCauchy} do not converge, which is the same as saying that there is no real number $c_0$ for which \eqref{NDA} holds true. This follows from Theorem 1 in \cite{GabettaRegazzini2012}. With regard to this introductory discussion, it is worth recalling the following

\begin{proposition}\label{prop1}
In order that the solution of \eqref{eqCauchy} converge weakly to a p.m. on $\R$, it is necessary that
\begin{equation}\label{rho}
\limsup_{x\to+\infty}x^\a(1-F^*_0(x))<+\infty.
\end{equation}  
\end{proposition}

This implies that convergence of the solution of \eqref{eqCauchy} is incompatible with unboundedness of the function 
\[
\rho(x):=x^\a(1-F^*_0(x))\qquad(x>0).
\]
This statement has been originally formulated in the first part of Step 2 in Section 3 of \cite{GabettaRegazzini2012}, but the argument used therein is not complete. Hence, we seize the opportunity of the present discussion to provide, in  Appendix B, a new proof of Proposition \ref{prop1}. As a matter of fact, some progress in the discovery of suitable necessary conditions for the explosion is achieved in the following proposition, according to which the family $\{\mu(\cdot,t)\}_{t\geq0}$ is uniformly tight if $\limsup_{x\to+\infty}\rho(x)$ turns out to be finite. Expressed in a form which is of use to the present discussion, uniform tightness means that for every strictly positive $\veps$ there is a bounded, closed interval $I_\veps$ of $\R$ such that $\inf_{t\geq0}\mu(I_\veps,t)>1-\veps$.

\begin{proposition}\label{prop2}
If $\limsup_{x\to+\infty}x^\a(1-F^*_0(x))<+\infty$, then $\{\mu(\cdot,t)\}_{t\geq0}$ is uniformly tight.
\end{proposition}

See Appendix C for the proof.\newline

As a consequence, explosion can not occur if $\rho$ is bounded. The question arises whether the unboundedness of $\rho$ is also sufficient in order that the solution of \eqref{eqCauchy} explode. As things stand, we are unable to state whether this is true or false. On the other hand, with the help of a slightly stronger condition, we manage to state explosion and to quantify its rate. The condition at issue $-$ which is reminiscent of the \textit{naif} idea of ultra-heavy-tailed distribution mentioned in Section \ref{sec:intro} $-$ reads

\textit{There are a strictly increasing  and divergent sequence {\rm $(x_m)_{m\geq1}$} of positive real numbers such that}
\begin{equation}\label{cond1}
\lim_{m\to+\infty}x_m^\a[1-F_0^*(x_m)]=+\infty
\end{equation}
\textit{and a nonincreasing function $c\colon \R^+\to (0,\a]$ for which:} 
\[
\theta_m=c_m x_m^\a[1-F^*_0(x_m)]\stackrel{m\to+\infty}{\longrightarrow}+\infty
\]
\textit{for $c_m:=c(x_m)$ and, for every $m\geq1$,}
\begin{equation}\label{cond2}
\begin{split}
&x\longmapsto F_0^*(x)-\dfrac{\theta_m}{\a(-x)^\alpha}\textit{ is nondecreasing on $(-\infty,-x_m)$.}
\end{split}
\end{equation}
For a better understanding of \eqref{cond2}, it is worth reformulating it when $F^*_0$ is absolutely continuous, with probability density function $f^*_0$. In fact, in such a case, \eqref{cond2} turns out to be equivalent to 
\begin{equation}\label{cond2d}
x\mapsto f^*_0(x)\geq \frac{\theta_m}{|x|^{\a+1}}\quad\text{a.e. on $(-x_m,x_m)^c$}
\end{equation}
for every $m=1,2,\dots$. 

\begin{examples}\label{es1}
{\rm
The aim here is to provide a few significant initial p.d.'s, belonging to domains of attraction of stable laws or to domains of attraction of extreme value distributions, which satisfy conditions \eqref{cond1}-\eqref{cond2}.
\begin{itemize}
\item[(1)] For the probability density function 
\[
f^*_0(x)=\dfrac{\beta}{2}\dfrac{1}{|x|^{\beta+1}}\I_{\{|x|\geq1\}}
\]
with $\beta<\a$, condition \eqref{cond2d} becomes $\beta x^{\a-\beta}\geq c_m x_m^{\a-\beta}$ for every $x\geq x_m$ which is satisfied by taking $c(\cdot)$ identically equal to some strictly positive constant $c\leq\beta$.
\item[(2)] For any p.d.f. $F^*_0$ such that $F^*_0(x)=e^{-1/x^\beta}$ whenever $x\geq1$ and $\beta<\a$, condition \eqref{cond2d} reduces to 
\[
c_m\leq \dfrac{\beta x^{\a-\beta}e^{-1/x^\beta}}{x_m^\a(1-e^{-1/x_m^\beta})}\qquad (x\geq x_m)
\]
which is satisfied if $c_m\leq \beta x_m^{-\beta}/(e^{x_m^{-\beta}}-1)$ and, \textit{a fortiori}, if $c(\cdot)\equiv c\leq \beta(1-\delta)$ for some $0<\delta<1$. In turn, this last inequality holds when, without loss of generality, we suppose that $x_1$ is sufficiently large.
\item[(3)] The previous initial p.d.f.'s satisfy one of the so-called Von Mises conditions, that is, their probability density functions $f^*_0$ meet
\begin{equation}\label{ipotesidensita}
\lim_{x\to+\infty}\dfrac{x f^*_0(x)}{1-F^*_0(x)}=\beta<\a.
\end{equation}
When \eqref{ipotesidensita} is in force, condition \eqref{cond2d} becomes
\[
c_m\leq \dfrac{x f^*_0(x)}{1-F^*_0(x)}\cdot \dfrac{x^\a[1-F^*_0(x)]}{x_m^\a[1-F^*_0(x_m)]} \qquad (x\geq x_m).
\]
Since \eqref{ipotesidensita} implies that $x\mapsto x^\a[1-F^*_0(x)]$ is nondecreasing for $x\to+\infty$, then condition \eqref{cond2d} holds if $c(\cdot)$ reduces to a constant $c\leq \beta-\delta$ for some $0<\delta<\beta$ and $x_1$ is taken sufficiently large. 
\item[(4)] For any p.d.f. $F^*_0$ such that $F^*_0(x)=e^{-1/(\log x)^\beta}$ whenever $x\geq1$ and $\beta<\a$, condition \eqref{cond2d} reads 
\[
c_m\leq \beta \dfrac{e^{-1/(\log x)^\beta}x^\a}{(\log x)^{\beta+1}}\dfrac{1}{x_m^\beta(1-e^{-1/(\log x_m)^\beta})}\qquad (x\geq x_m)
\]
which is satisfied for sufficiently large values of $x_1$ and for $c(x)\leq \dfrac{\beta(1-\delta)}{\log x}$ whenever $x$ is greater than $1$ and $\delta$ is any constant in $(0,1)$.
\end{itemize}
}
\end{examples}

At this stage, we present the main result of the paper.

\begin{theorem}\label{teorema}
Let $\mu(\cdot,t)$ be the solution of \eqref{eqCauchy} with symmetrized initial p.d.f. $F^*_0$ satisfying \eqref{cond1}-\eqref{cond2}. Then 
\[
\lim_{t\rightarrow +\infty}\mu((-\infty,-x),t)=\lim_{t\rightarrow +\infty} \mu((x,+\infty),t)=\dfrac{1}{2}\quad(x\geq0).
\]
Moreover, 
\[
\lim_{t\rightarrow+\infty}\mu \left( \left[ \rho(x_{m(t)})^{1/\alpha}\veps(t),+\infty\right),t\right)=\dfrac{1}{2}
\]
where $\veps\colon\R^+\to\R^+$ is any positive function which vanishes at infinity and $m:\R^+\to\R^+$ is defined by 
\[
m(t):=\sup \left\lbrace n\in\mathbb{N} : x_n\leq 2^{1/2}e^{\tau' t}\wedge(H_1)^{-1}(e^{-\tau_1''t})\wedge(H_{2-\delta})^{-1}(e^{-\tau_2''t})\right\rbrace
\]
where $\tau'$, $\tau_1''$ and $\tau_2''$ are suitable strictly positive constants, $\delta$ is any constant in $(0,1)$ and $H_p(x):=c(x)[1-F^*_0(x)]^p$.

\end{theorem}

The constant $\tau'$ is equal to $\tau\sigma/2$, where $\tau$ satisfies $2-\a-2\tau>0$ and $\sigma$ is a suitable fixed positive number, that will be specified in the following brief description of the proof of the theorem. As for $\tau_1''$ and $\tau_2''$, they are specified at the end of Section \ref{sec:prove}. As usual, given any nonincreasing function $H$ on $\R$, $H^{-1}$ denotes its \textit{generalized inverse}, i.e., $H^{-1}(u):=\inf\{x:\; H(x)< u\}$, for every $u$ in the interior of the convex hull of the range of $H$. 

The core of the argument used, in Section \ref{sec:prove}, to prove the theorem, reduces to verifying the validity of the following bound, for every $x>0$,
\begin{equation}\label{megabound}
\begin{split}
2\mu([0,x],t)&\leq \dfrac{2}{\pi}\Gamma\left( \dfrac{1}{\alpha}\right)\Big[\dfrac{(1-F_0^*(x_{m(t)}))^{\frac{1-\delta}{\alpha}}}{\alpha(1+2\pi)^2}+\dfrac{x}{a_{m(t)}^{1/\alpha}}\Big]+ C^*_1 (1-F_0^*(x_{m(t)}))^{\delta}\\
&+C^*_{2} \dfrac{x_m^{3\a}e^{-(a_{m(t)}-\eta)x_{m(t)}^\alpha}}{(a_{m(t)}-\eta)}+ C^*_3 \dfrac{e^{-(2-\alpha-2\tau+\alpha \tau/2)\sigma t}}{a_{m(t)}-\eta}+ C^*_4 e^{-\gamma t}
\end{split}
\end{equation}
where $C^*_1,\dots,C^*_4$ are positive constants, and
\begin{equation}\label{am}
a_m:=\frac{2}{\a} \dfrac{\a \theta_m x_m^\a}{2\theta_m+\a x_m^\a[2F^*_0(x_m)-1]}\int_0^{+\infty} \dfrac{\sin x}{x^\alpha}dx\quad(m=0,1,\dots)
\end{equation}
\[
\gamma:=\min\Big\{1-2R_q-\dfrac{\tau_2''q}{2},1-2R_4-2\tau_1'',\Lambda-\dfrac{\tau_1''}{\a},1-2R_q-\dfrac{q\sigma\a}{2}-2\tau_1''\Big\}
\]
\[
\text{$q$ is any fixed number in $(2,+\infty)$ and  $R_q:=\frac{\Gamma\left( \frac{q}{2}+\frac{1}{2}\right) }{\sqrt{\pi}\Gamma\left( \frac{q}{2}+1\right)}$}
\] 
\[
\text{$\sigma$ is any number in $(0, \frac{2(1-2R_q)}{q\a})$}
\]
\[
\Lambda:=\min\Big\{\sigma, 1-2R_q-\dfrac{q\sigma\a}{2}\Big\}. 
\]

\begin{remarks}\label{oss1}
{\rm 
\begin{itemize}
\item[(a)] With a view to the understanding of the explosion and the quantification of its rate, one notes that $t\mapsto m(t)$ diverges monotonically as $t$ goes to infinity and so does $t\mapsto \rho(x_{m(t)})$ by taking, if necessary, a subsequence of $(x_m)_{m\geq1}$ in the definition of $t\mapsto m(t)$. A nice fact is the transparent and simple connection of the lower bound to the rate of explosion, given by $\rho(x_{m(t)})^{1/\alpha}\veps(t)$, with the "speed" of divergence to infinity of the LHS of \eqref{cond1}. This speed, in its turn, can be viewed as an index of the departure of data, which cause explosion, from those that $-$ according to Proposition \ref{prop2} $-$ entail \textit{relative compactness} of the family $\{\mu(\cdot,t):\;t\geq0\}$ associated with the solution of \eqref{eqCauchy}. Examples \ref{es2} will illustrate this kind of ideas. See Section 5 of \cite{Billingsley} for the definition of the above relative compactness. In view of the arbitrariness of $\veps(\cdot)$, it is of course desirable to choose forms which decay to $0$ as slowly as possible.
\item[(b)] It is easy to check that if $F^*_0$ is absolutely continuous and \eqref{cond1},\eqref{cond2d} hold, then $\lim_{x\to+\infty}x^\a[1-F^*_0(x)]=+\infty$.
\item[(c)] When $\lim_{x\to+\infty}x^\a[1-F^*_0(x)]=+\infty$, we may and do replace \eqref{cond2} with the condition
\[
\begin{split}
&x\longmapsto F_0^*(x)-\dfrac{\theta(u)}{\a(-x)^\alpha}\textit{ is nondecreasing on $(-\infty,-u)$ for every $u>0$}
\end{split}
\]
where $\theta(u):=c(u)u^\a[1-F^*_0(u)]$. Then, the main portion of Theorem \ref{teorema} can be reformulated as:
\[
\lim_{t\to+\infty}\mu\Big([\rho(u(t))^{1/\a}\veps(t),+\infty),t\Big)=\dfrac{1}{2}
\]
where $u(t):=\min\{2^{1/2}e^{\tau' t}\wedge(H_1)^{-1}(e^{-\tau_1''t})\wedge(H_{2-\delta})^{-1}(e^{-\tau_2''t})\}.$ Moreover, if a constant version for $c(\cdot)$ is admissible, then 
\begin{equation}\label{boundU}
u(t)\geq C e^{tA}
\end{equation} 
for some positive constant $C$ and $A=\tau'\wedge\tau_1''/\a\wedge\tau''_2/\a(2-\delta)$. The last inequality is valid since, if $c(\cdot)$ is a constant function, then 
\[
u(t)=\min\{\sqrt{2}e^{t\tau'},(F^*_0)^{-1}(1-e^{-t\tau''})\}
\]
with $\tau'':=\min\{\tau''_1,\tau''_2/(2-\delta)\}$. Moreover, since $1-F^*_0(x)\geq 1/x^\a$, for $x\to+\infty$, we get $F^*_0(1-e^{-t\tau''})\geq e^{t\tau''/\a}$ for $t\to+\infty$.
\end{itemize}
}
\end{remarks}

\begin{examples}\label{es2}
{\rm Here we use Theorem \ref{teorema} to estimate the rate of explosion of the solution of equation \eqref{eqCauchy} in a couple of cases.
\begin{itemize}
\item[(1)] If $F^*_0$ is the same as in Examples \ref{es1}.1, then $(F^*_0)^{-1}(y)=\dfrac{1}{(2(1-y))^{1/\beta}}$ for $y\to 1^-$. Hence, $u(t)=Ce^{tB}$ with $B=\tau'\wedge \tau''/\beta$, as $t\to+\infty$, which entails that
\[
x(t)=\rho(u(t))^{1/\a}\veps(t)=C' \veps(t) e^{tB(\a-\beta)/\a}
\]
for some positive constant $C'$, as $t\to+\infty$.
\item[(2)] Next, consider any $F^*_0$ that, for $x\to+\infty$, is equal to $F^*_0(x)=1-\log x/x^\a$. It is easy to check that this initial p.d.f. satisfies \eqref{cond1}-\eqref{cond2}. Since $\log x\leq x^\veps$ as $x\to+\infty$, for any $\veps>0$, it is easy to check that $(F^*_0)^{-1}(1-e^{-t\tau''})\leq e^{t\tau''/(\a-\veps)}$, as $t\to+\infty$, and hence there exists a constant $c>0$ such that $u(t)\leq c e^{tD}$ with $D:=\tau'\wedge \tau''/(\a-\veps)$, as $t\to+\infty$. Combining this bound with \eqref{boundU}, we get
\[
c' t^{1/\a}\veps(t)\leq x(t)\leq c'' t^{1/\a}\veps(t)
\]
for suitable constants $0<c'<c''$, as $t\to+\infty$.
\end{itemize}

}
\end{examples}

These examples show that both exponential and nonexponential lower bounds for the rate of explosion may occur, depending on the specific initial data. Furthermore, from Appendix D, it turns out that the rates of explosion $a_t$, in these very same examples, satisfy the following inequalities, respectively,
\[
C' \veps(t) e^{tB(\a-\beta)/\a}\leq a_t< \dfrac{C'}{\veps_1(t)^{1/\beta}}e^{t[2R_{2\beta/\a}-1+B(\a-\beta)\beta\a]/\beta}
\]
\[
c't^{1/\a}\veps(t)\leq a_t< c''\frac{t^{1/\a}}{\veps_1(t)}
\]
where $t\mapsto \veps_1(t)$ is any strictly positive function on $(0,+\infty)$ such that $\veps_1(t)\searrow 0$, as $t\to+\infty$.

\section{Proof of Theorem \ref{teorema}}\label{sec:prove}
The proof is split into steps, the first of which is concerned with an enlargement of the probability space $(\O,\CF,\CP_t)$ to support some new random elements. For such an enlargement, we preserve the previous symbol without possibility of ambiguity. \newline
\textbf{Step 1: }\textit{Enlargement of the probabilistic framework.} The new elements, to be considered together with $\tte$, $\ti$, $\tnu$, are denoted by 
\begin{itemize}
\item $\tb=(\tb_n)_{n\geq1}$
\item $U=(U_n)_{n\geq1}$
\item $Z=(Z_n)_{n\geq1}$
\end{itemize} 
and they are defined in combination with a family $\{\ptm:\; t>0,\; m=1,2,\dots\}$ of p.m.'s on $(\O,\CF)$ in such a way they meet the following conditions:
\begin{itemize}
\item $\tte$, $\ti$, $\tnu$, $\tb$, $U$, $Z$ are mutually independent with respect to each $\ptm$.
\item $\tte$, $\ti$, $\tnu$ preserve, under each $\ptm$, the same distributions as those defined at the beginning of Section \ref{sec:preliminaries}.
\item $b$ is a Bernoulli sequence with parameter $1/K_{1,m}$, with respect to each $\ptm$, where $1/K_{1,m}:=2\theta_m/(\a x_m^\a)+2F^*_0(x_m)-1$.
\item both the $U$'s and the $Z$'s are, under each $\ptm$, i.i.d. random numbers with common p.d.f.'s $G_{1,m}$ and $G_{2,m}$, respectively, defined by\newline
\[
\begin{split}
&G_{1,m}(x):=\I_{(-\infty,-x_m)}(x)\Big[\dfrac{K_{1,m}\theta_m}{\a}\cdot \frac{1}{(-x)^\alpha} \Big]+\I_{[-x_m,x_m)}(x)K_{1,m}\Big[F_0^*(x)\\
&\qquad\qquad\qquad+\dfrac{\theta_m}{\a x_m^\a}+F^*_0(x_m)-1\Big]+\I_{[x_m,+\infty)}(x)\Big[1-\dfrac{K_{1,m}\theta_m}{\a}\cdot \frac{1}{x^\alpha}\Big]
\end{split}
\]
and
\[
\begin{split}
G_{2,m}(x)&:=\I_{(-\infty,-x_m)}(x)K_{2,m}\Big[F_0^*(x)- \dfrac{\theta_m}{\a(-x)^\a}\Big]+\frac{1}{2}\I_{[-x_m,x_m)}(x)\\
&+\I_{[x_m,+\infty)}(x)K_{2,m}\Big[F_0^*(x)+\dfrac{\theta_m}{\a x^\a}+1-2F_0^*(x_m)-\frac{2\theta_m}{\a x_m^\a}\Big].
\end{split}
\]
\end{itemize} 
To verify that they are p.d.f.'s, it is enough to recall, apropos of $G_{2,m}$, condition \eqref{cond2}. Moreover, the random numbers $\tb_jU_j+(1-\tb_j)Z_j$ ($j=1,2,\dots$) turn out to be i.i.d. with common p.d.f. $F^*_0$, with respect to each $\ptm$. Hence, the c.f. of the random sum $\sum_{j=1}^{\tnu}|\tbe_{j,\tnu}|(\tb_jU_j+(1-\tb_j)Z_j)$, under each $\ptm$, is the solution of problem \eqref{eqCauchy} with initial datum $\vphi^*_0$. Thus, we can now aim at verifying the validity of \eqref{megabound} under the additional information that the LHS therein coincides with
\begin{equation}\label{bound}
\CP_{t,m(t)}\Big\{-x\leq \sum_{j=1}^{\tnu}|\tbe_{j,\tnu}|(\tb_jU_j+(1-\tb_j)Z_j)\leq x\Big\}.
\end{equation}
\newline
\textbf{Step 2: }\textit{Analysis of \eqref{bound}.} We begin by verifying the following two inequalities, in which the symbol comma between two conditions stands for their intersection: 
\begin{equation*}
\begin{split}
&\ptm\Big\{-x\leq \sum_{j=1}^{\tnu}|\tbe_{j,\tnu}|(\tb_jU_j+(1-\tb_j)Z_j)\leq x\Big\}\\
&\leq \ptm\Big\{-x\leq \sum_{j=1}^{\tnu}|\tbe_{j,\tnu}|U_j+ \sum_{j=1}^{\tnu}|\tbe_{j,\tnu}|(1-\tb_j)Z_j-\sum_{j=1}^{\tnu}|\tbe_{j,\tnu}|(1-\tb_j) U_j\leq x,\\ 
&\qquad\qquad\qquad\qquad\qquad\qquad\qquad\qquad\qquad\qquad\Big|\sum_{j=1}^{\tnu}|\tbe_{j,\tnu}|(1-\tb_j) U_j\Big|\leq A_m\Big\}
\end{split}
\end{equation*}
\begin{equation}\label{2termini}
\begin{split}
&+ \ptm\Big\{\Big|\sum_{j=1}^{\tnu}|\tbe_{j,\tnu}|(1-\tb_j) U_j\Big|> A_m\Big\}\\
&\leq \ptm\Big\{\Big|\sum_{j=1}^{\tnu}|\tbe_{j,\tnu}|(1-\tb_j) U_j\Big|> A_m\Big\}\\
&+ \ptm\Big\{-x-A_m\leq \sum_{j=1}^{\tnu}|\tbe_{j,\tnu}|U_j+ \sum_{j=1}^{\tnu}|\tbe_{j,\tnu}|(1-\tb_j)Z_j\leq x+A_m\Big\}\\
& =: \CL_1+\CL_2.
\end{split}
\end{equation}
Throughout the paper, we preserve the symbols $\theta_m:=c_m x_m^\a[(1-F^*_0(x_m))]$, $a_m:=\dfrac{2}{\a}K_{1,m}\theta_m\int_0^{+\infty} (\sin x)/x^\alpha dx$, $k_m:=a_m^{1/\alpha}(1-F_0^*(x_m))^{\frac{1-\delta}{\alpha}}$ for $\delta$ in $(0,1)$, $A_m:=\frac{k_m}{(1+2\pi)^2}$, $m\geq1$. \newline
\newline

\textit{Bound for $\CL_1$.} We prove that 
\begin{equation}\label{tesidoob}
\CL_1\leq \dfrac{1}{\a+1}\frac{k_m(a_m+\eta)}{A_m^{\a+1}K_{2,m}}+\ptm\{\max_{j=1,\dots,\tnu}|\tbe_{j,\tnu}(1-\tb_j)|>\veps\}
\end{equation}
holds for every $\veps$ not greater than a suitable positive $\veps_m$ and $\eta$ in $(0,a_m)$. Given the $\sigma$-algebra $\CS$ generated by $(\tnu,\tte,\ti,\tb)$, we write $\CL_1=\E_{t,m}[\CA_{1,\veps}+\CA_{2,\veps}]$ where 
\[
\begin{split}
\CA_{1,\veps}&=\CA_{1,\veps}(\o)\\
&:=\ptm\Big\{\Big|\sum_{j=1}^{\tnu}|\tbe_{j,\tnu}|(1-\tb_j) U_j\Big|> A_m \Big|\CS\Big\}(\o)\I_{\{\max_{j}|\tbe_{j,\tnu}(1-\tb_j)|\leq \veps\}}(\o)
\end{split}
\]
and
\[
\begin{split}
\CA_{2,\veps}&=\CA_{2,\veps}(\o)\\
&:=\ptm\Big\{\Big|\sum_{j=1}^{\tnu}|\tbe_{j,\tnu}|(1-\tb_j) U_j\Big|> A_m\Big|\CS\Big\}(\o)\I_{\{\max_{j}|\tbe_{j,\tnu}(1-\tb_j)|> \veps\}}(\o)
\end{split}
\] 
for every $\o$ in $\O$. We proceed to bound $\CA_{1,\veps}$ by using the inequality 
\begin{equation}\label{doob}
\ptm\{|Y|>C\}\leq\frac{1}{\Delta}\Big(1+\frac{2\pi}{C\Delta}\Big)^2\int_0^\Delta\Big[1-\Re \vphi_Y(u)\Big]du
\end{equation}
where $Y$ is any random number and $\vphi_Y$ its c.f.. For a proof of this inequality see Subsection 8.3 of \cite{ChowTeicher}. Then, putting $C=\Delta^{-1}=A_m$, we obtain
\begin{equation}\label{passodoob}
\begin{split}
&\ptm\Big\{\Big|\sum_{j=1}^{\tnu}|\tbe_{j,\tnu}|(1-\tb_j) U_j\Big|> A_m\Big|\CS\Big\}\\
&\qquad\qquad\qquad\qquad\qquad\leq k_m\int_{( 0,A_m^{-1}) } \Big( 1-\prod_{j=1}^{\tnu} \hat{g}_{1,m}( |\tbe_{j,\tnu}|(1-\tb_j)\xi)\Big) d\xi
\end{split}
\end{equation}
where $\hat{g}_{1,m}(\xi):=\int_\R e^{i\xi x}G_{1,m}(dx)$ is a real-valued function in view of the symmetry of $G_{1,m}$. Now we give a bound for the RHS of \eqref{passodoob} based on the following lemma, whose proof is deferred to Appendix E.

\begin{lemma}\label{g1m}
The p.d.f. $G_{1,m}$ belongs to the standard domain of attraction of the stable law with c.f.
\[
\xi\mapsto\exp(-a_m|\xi|^\alpha)\quad (\xi\in\R),
\]
$a_m$ being the same as in \eqref{am}. Moreover, the following equality
\begin{equation*}
1-\hat{g}_{1,m}(\xi)=(a_m+v_m(\xi))\xi^\alpha
\end{equation*}
holds for every strictly positive $\xi$, with $v_m$ defined by
\[
\begin{split}
v_m(\xi)=&\dfrac{2}{\xi^\alpha}\int_0^{\xi x_m} \sin x \Big[1-K_{1,m}\Big(F^*_0\Big(\dfrac{x}{\xi}\Big)+F^*_0(x_m)-1+\dfrac{\theta_m}{\a x_m^\a} \Big)\Big]dx\\
&-\dfrac{2}{\a}K_{1,m}\theta_m\int_0^{\xi x_m}\dfrac{\sin x}{x^\alpha}dx\quad(\xi>0).
\end{split}
\]
For each $m$, $v_m$ is $o(1)$ as $\xi\to 0^+$ and its supremum norm $\left \|v_m\right \|$ satisfies
\[
\left \|v_m\right \|\leq 2^{1-\a}x_m^\a +\frac{2}{\a}K_{1,m}\theta_m\int_{(0,\pi)}\frac{\sin x}{x^\a}dx.
\]
\end{lemma}
Coming back to \eqref{passodoob}, setting $q_j:=|\tbe_{j,\tnu}|(1-\tb_j)$, $j=1,\dots,\tnu$, we write 
\begin{equation}\label{fi}
\begin{split}
\varphi_{\tnu,1,m}(\xi)&:=\prod_{j=1}^{\tnu} \hat{g}_{1,m}(q_{j} \xi)\\
& = \exp \left[ \sum_{j=1}^{\tnu} \log(1-(1- \hat{g}_{1,m}(q_{j} \xi)))\right]\\
&= \exp \left[ -\sum_{j=1}^{\tnu}\left( 1-\hat{g}_{1,m}(q_{j} \xi)-\dfrac{4}{5}\theta_{j}(1-\hat{g}_{1,m}(q_{j} \xi))^2 \right) \right]
\end{split}
\end{equation}
which, in view of Lemma 3 of Section 9.1 in \cite{ChowTeicher}, is valid for every $\xi$ in a suitable neighbourhood of the origin such that $0\leq1-\hat{g}_{1,m}(q_{j} \xi)\leq 3/8$ for every  $j=1,\dots,\tnu$. A possible choice of the this neighbourhood can be derived from Lemma \ref{g1m}, that is 
\begin{equation}\label{2}
1-\hat{g}_{1,m}(q_{j} \xi)=(a_m+v_m(q_j\xi))|\xi|^\alpha q_j^\alpha \leq M_m q_{j}^\alpha|\xi|^\alpha 
\end{equation}
with  $M_m:=a_m+\left \|v_m\right \|$. Thus, the neighbourhood $\CN_d:=\{\xi\in\R:\;|\xi|^\a q_{(\tnu)}^\a\leq  3d/(8M_m)\}$, with $q_{(\tnu)}:=\max\lbrace q_{1},...,q_{\tnu}\rbrace$ and $d$ any number in $(0,1)$, serves our purpose. Hence, 
\[
\begin{split}
\varphi_{\tnu,1,m}(\xi)&=\exp \left[ -\sum_{j=1}^{\tnu}\left( 1-\hat{g}_{1,m}(q_{j} \xi)-\dfrac{4}{5}\theta_{j}(1-\hat{g}_{1,m}(q_{j} \xi))^2 \right) \right]\\
&=\exp \left[ -\sum_{j=1}^{\tnu} (a_m+v_m(\xi q_{j}))|\xi|^\alpha q_{j}^\alpha +\dfrac{4}{5}\sum_{j=1}^{\tnu}\theta_{j}(1-\hat{g}_{1,m}(q_{j} \xi))^2 \right]\\
&= \exp \left[ -a_m|\xi|^\alpha \sum_{j=1}^{\tnu}  q_{j}^\alpha\right] \exp \left[ -B_{\tnu}(\xi)+R_{\tnu}(\xi)\right]
\end{split}
\]
where $B_{\tnu}(\xi):=|\xi|^\alpha \sum_{j=1}^{\tnu}q_{j}^\alpha v_m(\xi q_{j})$, $R_{\tnu}(\xi):=\dfrac{4}{5}\sum_{j=1}^{\tnu}\theta_{j}(1-\hat{g}_{1,m}(q_{j} \xi))^2$ and $|\theta_j|\leq1$ for every $j$ and $\xi$ in $\CN_d$. Now, for every $\xi$ in $\CN_d$, \eqref{2} yields
\[
|R_{\tnu}(\xi)|\leq \dfrac{4}{5}M_m^2|\xi|^{2\alpha}\sum_{j=1}^{\tnu}  q_{j}^{2\alpha}
\]
and, setting $\bar{v}_m (\xi):=\sup_{0\leq x\leq\xi} |v_m(x)|$, 
\[
\begin{split}
\left| B_{\tnu}(\xi)\right| +\left| R_{\tnu}(\xi)\right|& \leq |\xi|^\alpha \sum_{j=1}^{\tnu} \bar{v}_m(\xi q_{(\tnu)})q_{j}^\alpha +\dfrac{4}{5}M_m^2|\xi|^{2\alpha}q_{(\tnu)}^\alpha \sum_{j=1}^{\tnu} q_{j}^\alpha \\
&\leq |\xi|^\alpha \sum_{j=1}^{\tnu} q_{j}^\alpha\left( \bar{v}_m(\xi q_{(\tnu)})+\dfrac{4}{5}M_m^2|\xi|^{\alpha}q_{(\tnu)}^\alpha \right).
\end{split}
\]
Since $|v_m|$ is $o(1)$ as $\xi\to 0^+$, for every $\eta$ in $(0,a_m)$, $d=d^*$ can be chosen sufficiently small so that $\bar{v}_m(\xi q_{(\tnu)})+\dfrac{4}{5}M_m^2|\xi|^{\alpha}q_{(\tnu)}^\alpha\leq \eta$ holds true for every $\xi$ in $\CN_{d^*}$ and, therein,
\begin{equation}\label{somma}
\left| B_{\tnu}(\xi)\right| +\left| R_{\tnu}(\xi)\right|\leq \eta |\xi|^\alpha \sum_{j=1}^{\tnu} q_{j}^\alpha.
\end{equation}
Now, for any $\veps\leq\veps_m:=A_m\Big(\dfrac{3d^*}{8M_m}\Big)^{1/\a}$, \eqref{passodoob} and \eqref{fi} give
\[
\begin{split}
\CA_{1,\veps}&\leq \I_{\{q_{(\tnu)}\leq\veps\}}k_m\int_{\left( 0,A_m^{-1}\right) }  \left( 1-\exp {\left[ -a_m |\xi|^\alpha \sum_{j=1}^{\tnu}  q_{j}^\alpha -B_{\tnu}(\xi)+R_{\tnu}(\xi)\right]}\right)d\xi\\
& \leq \I_{\{q_{(\tnu)}\leq\veps\}}k_m\int_{\left( 0,A_m^{-1}\right) } \left(a_m|\xi|^\alpha \sum_{j=1}^{\tnu}  q_{j}^\alpha +B_{\tnu}(\xi)-R_{\tnu}(\xi)\right)d\xi\\
&\qquad\qquad\qquad\qquad\qquad\qquad\qquad\qquad\qquad\qquad\text{(from $1-e^{-x}\leq x$)}\\
& \leq \I_{\{q_{(\tnu)}\leq\veps\}}k_m\int_{\left( 0,A_m^{-1}\right) } \left(a_m|\xi|^\alpha \sum_{j=1}^{\tnu}  q_{j}^\alpha +\eta |\xi|^\alpha \sum_{j=1}^{\tnu} q_{j}^\alpha\right)d\xi\\
&\qquad\qquad\qquad\qquad\qquad\qquad\qquad\qquad\qquad\qquad\text{(from \eqref{somma}, since $\veps\leq\veps_m$)}\\
&= \I_{\{q_{(\tnu)}\leq\veps\}}k_m\frac{(a_m+\eta)}{\a+1}\dfrac{1}{A_m^{\a+1}}\sum_{j=1}^{\tnu} q_{j}^\alpha.
\end{split}
\]
Then $\CL_1\leq \CE_1+\CE_2$ where $\CE_1 := k_m\frac{(a_m+\eta)}{\a+1}\dfrac{1}{A_m^{\a+1}}\E_{t,m}\Big(\sum_{j=1}^{\tnu} q_{j}^\alpha\Big) $, $\CE_2:= \ptm\{q_{(\tnu)}>\veps\}$. As for $\CE_1$,
\[
\E_{t,m}\Big(\sum_{j=1}^{\tnu} q_{j}^\alpha\Big)=\sum_{j=1}^{\tnu} \E_{t,m}\left[ |\tbe_{j,\tnu}|^\alpha\right] \E_{t,m} (1-\tb_j)=\dfrac{1}{K_{2,m}}
\]
where the last equality is obtained recalling that $\sum_{j=1}^{\tnu} |\tbe_{j,\tnu}|^\alpha=1$ almost surely. Whence,
\[
\CE_1= \dfrac{1}{\a+1}\frac{k_m(a_m+\eta)}{A_m^{\a+1}K_{2,m}}
\]
which completes the proof of \eqref{tesidoob}.\newline

\textit{Bound for $\CL_2$.} Let $F_{1,m}$ and $F_{2,m}$ be conditional p.d.f.'s for $\sum_{j=1}^{\tnu} |\tbe_{j,\tnu}|U_j$ and $\sum_{j=1}^{\tnu} |\tbe_{j,\tnu}|(1-\tb_j)Z_j$ respectively, given $\CS$. Moreover, let $S_{a_m}$ be the p.d.f. associated with the stable law mentioned in Lemma \ref{g1m}. Then,
\begin{equation}\label{L2}
\begin{split}
\CL_2&= \E_{t,m}\left[ \int_\mathbb{R} F_{1,m}\left(x+A_m-s\right)-F_{1,m}\left(-x-A_m-s\right)dF_{2,m}(s)\right]\\
& \leq\E_{t,m}\left[ \int_\mathbb{R}\left| F_{1,m}\left(x+A_m-s\right)-S_{a_m}\left(x+A_m-s\right)\right|dF_{2,m}(s)\right]\\
&+ \E_{t,m}\left[ \int_\mathbb{R}\left| S_{a_m}\left(x+A_m-s\right)-S_{a_m}\left(-x-A_m-s\right)\right|dF_{2,m}(s)\right]\\
&+ \E_{t,m}\left[ \int_\mathbb{R}\left| F_{1,m}\left(-x-A_m-s\right)-S_{a_m}\left(-x-A_m-s\right)\right|dF_{2,m}(s)\right]\\
&\qquad\\
& =: \CL_{2,1}+\CL_{2,2}+\CL_{2,3}.
\end{split}
\end{equation}
We start with $\CL_{2,2}$ by recalling that $S_{a_m}$ is symmetric and unimodal (cf. Theorem 2.5.3 in \cite{IbragimovLinnik}). The integrand function represents a probability on an interval of length $2( x+A_m)$ so that, keeping in mind the concept of \textit{L\'evy's concentration function}, 
\[
\begin{split}
&\left| S_{a_m}\left( x+A_m-s\right) -S_{a_m}\left(-x-A_m-s\right)\right|\\
&\qquad\qquad\qquad \leq S_{a_m}\left( x+A_m\right) -S_{a_m}\left(-x-A_m\right)\\
&\qquad\qquad\qquad =\dfrac{1}{\pi} \int_\mathbb{R} \dfrac{1}{t}\sin\left( t\left( x+A_m\right) \right) e^{-a_m|t|^\alpha}dt\\
&\qquad\qquad\qquad\qquad\qquad\qquad\qquad\qquad\qquad\quad \text{(from L\'evy's inversion formula)}\\
&\qquad\qquad\qquad \leq \dfrac{2}{\pi}\left( x+A_m\right)\int_{(0,+\infty)} e^{-a_m|t|^\alpha}dt\\
&\qquad\qquad\qquad =\dfrac{2}{\pi}\left( x+A_m\right)\dfrac{1}{\alpha a_m^{1/\alpha}}\Gamma\left( \dfrac{1}{\alpha}\right).
\end{split}
\]
Whence,
\begin{equation}\label{secondo}
\CL_{2,2}\leq \Gamma\left( \dfrac{1}{\alpha}\right)\dfrac{2}{\a\pi}\dfrac{ x+A_m}{ a_m^{1/\alpha}}.
\end{equation}
We proceed to study $(\CL_{2,1}+\CL_{2,3})$. It is easy to check that 
\[
\CL_{2,1}+\CL_{2,3}\leq 2\E_{t,m}\left[ K(F_{1,m},S_{a_m})\right]
\]
where $K(F,G)$ stands for the \textit{Kolmogorov distance} between p.d.f.'s $F$ and $G$, i.e.
\[
K(F,G):=\sup_{x\in\R}|F(x)-G(x)|.
\] 
Since $F_{1,m}$ is p.d.f. of a weighted sum of i.i.d. random numbers whose common p.d.f. belongs to the standard domain of attraction of $S_{a_m}$, Proposition 3.4 in \cite{BaLaRe} can be applied to obtain
\begin{equation}\label{kolm}
\begin{split}
&K(F_{1,m},S_{a_m})\\
&\qquad\leq \dfrac{2}{\pi}\sum_{j=1}^{\tnu} |\tbe_{j,\tnu}|^\alpha \int_0^{\tilde{d}_m/\tbe_{(\tnu)}} \exp{\left( -(a_m-\eta)\xi^\alpha\right)} \xi^{\alpha-1} H_m(\xi,|\tbe_{j,\tnu}|)d\xi\\
&\qquad+\dfrac{8}{5\pi}M_m^2 N_{2\alpha,m}\sum_{j=1}^{\tnu} |\tbe_{j,\tnu}|^{2\alpha}+\dfrac{64}{25\pi}M_m^4 N_{4\alpha,m}\sum_{j=1}^{\tnu} |\tbe_{j,\tnu}|^{3\alpha}+\dfrac{\textbf{c}}{\tilde{d}_m}||g_m||\tbe_{(\tnu)}\\
&\qquad\\
&\qquad\Big(=: \CK_1+\CK_2+\CK_3+\CK_4\Big)
\end{split}
\end{equation}
where:
\[
{H}_m(\xi,s):=|v_m(\xi s)|(1+2|\xi|^\alpha |v_m(\xi s)|)
\]
\[N_{l,m}:=\int_0^{+\infty}\exp{(-(a_m-\eta)\xi^\alpha)}\xi^{l-1}d\xi\quad\text{and}\quad\tilde{d}_m=\left( \dfrac{3d^*}{8M_m}\right)^{1/\alpha};
\]
$g_m$ is a probability density function of $S_{a_m}$,  $\tbe_{(\tnu)}:=\max\{|\tbe_{1,\tnu}|,\dots,|\tbe_{\tnu,\tnu}|\}$ and $\textbf{c}= 24/\pi$ according to the formulation given in \cite{ChowTeicher} of a classical inequality due to Berry and Esseen.  As to $\CK_1$,
\[
\begin{split}
&\E_{t,m}(\CK_1)\\
&\quad\leq \dfrac{2}{\pi}\E_{t,m}\left[\sum_{j=1}^{\tnu} |\tbe_{j,\tnu}|^\alpha \int_0^{+\infty} \exp{\left( -(a_m-\eta)\xi^\alpha\right)} \xi^{\alpha-1} \sup_{0\leq s\leq \tbe_{(\tnu)}} H_m(\xi,s)d\xi\right]\\
&\quad =\dfrac{2}{\pi}\E_{t,m}\left[ \int_0^{+\infty} \exp{\left( -(a_m-\eta)\xi^\alpha\right)} \xi^{\alpha-1} \sup_{0\leq s\leq \tbe_{(\tnu)}} H_m(\xi, s)\I_{\lbrace\tbe_{(\tnu)}\leq e^{-\sigma t}\rbrace}d\xi\right]\\
&\quad+ \dfrac{2}{\pi}\E_{t,m}\left[ \int_0^{+\infty} \exp{\left( -(a_m-\eta)\xi^\alpha\right)} \xi^{\alpha-1} \sup_{0\leq s\leq \tbe_{(\tnu)}} H_m(\xi,s)\I_{\lbrace\tbe_{(\tnu)}> e^{-\sigma t}\rbrace}d\xi\right]\\
&\;\quad\qquad\qquad\qquad\qquad\qquad\qquad\qquad\qquad\qquad\qquad\qquad\qquad\qquad\qquad(\sigma>0)
\end{split}
\]
\[
\begin{split}
&\quad\leq\dfrac{2}{\pi}\int_0^{+\infty} \exp{\left( -(a_m-\eta)\xi^\alpha\right)} \xi^{\alpha-1} \sup_{0\leq s\leq e^{-\sigma t}} H_m(\xi,s)d\xi\\
&\quad + \dfrac{2}{\pi}\ptm(\tbe_{(\tnu)}> e^{-\sigma t})\int_0^{+\infty} \exp{\left( -(a_m-\eta)\xi^\alpha\right)} \xi^{\alpha-1} \sup_{s\geq 0} H_m(\xi,s)d\xi\\
&\quad\leq\dfrac{2}{\pi}\int_0^{+\infty} \exp{\left( -(a_m-\eta)\xi^\alpha\right)} \xi^{\alpha-1} \sup_{0\leq s\leq e^{-\sigma t}} H_m(\xi,s)d\xi\\
&\quad+ \dfrac{2}{\pi}e^{-t(1-q\sigma\alpha/2-2R_q)}\int_0^{+\infty} \exp{\left( -(a_m-\eta)\xi^\alpha\right)} \xi^{\alpha-1} \sup_{s\geq 0} H_m(\xi,s)d\xi
\end{split}
\]
where the last inequality follows from Lemma 1 in \cite{GabettaRegazzini2008} for every $q>0$ and $R_q:=\frac{\Gamma\left( \frac{q}{2}+\frac{1}{2}\right) }{\sqrt{\pi}\Gamma\left( \frac{q}{2}+1\right) }$. This, in combination with the definitions of $H_m(\cdot,\cdot\cdot)$ and $N_{l,m}$, gives 
\begin{equation}\label{k1}
\begin{split}
\E_{t,m}(\CK_1)& \leq\dfrac{2}{\pi}\int_0^{+\infty} \exp{\left( -(a_m-\eta)\xi^\alpha\right)} \xi^{\alpha-1}  \bar{H}_m(\xi,e^{-\sigma t})d\xi\\
& + \dfrac{2}{\pi}\left\|v_m\right\| \left(N_{\alpha,m} +2N_{2\alpha,m}\left\|v_m\right\|\right)e^{-t(1-q\sigma\alpha/2-2R_q)}
\end{split}
\end{equation}
where $\bar{H}_m(\xi,u):=\sup_{0\leq s\leq u}H_m(\xi,s)$. Now, recalling Proposition 8 in \cite{GabReg2006},
\begin{equation}\label{k23}
\E_{t,m}(\CK_2+\CK_3)\leq \dfrac{8}{5\pi}M_m^2N_{2\alpha,m}e^{-t(1-2R_4)}+\dfrac{64}{25\pi}M_m^4N_{4\alpha,m} e^{-t(1-2R_6)}.
\end{equation}
To bound $\E_{t,m}(\CK_4)$, we resort to Lemma 1 in \cite{GabettaRegazzini2008} to write
\begin{equation}\label{k4}
\E_{t,m}(\CK_4)\leq \dfrac{\textbf{c}\left \|g_{m}\right \|}{\tilde{d}_m}\left( e^{-\sigma t}+e^{-t(1-q\sigma\alpha /2-2R_q)}\right).
\end{equation}

At this stage, we combine \eqref{2termini} with \eqref{L2} and use \eqref{tesidoob} to bound $\CL_1$, \eqref{secondo} to bound $\CL_{2,2}$ and \eqref{k1}-\eqref{k4} to bound $\CL_{2,1}+\CL_{2,3}\leq \sum_{i=1}^4 \E_{t,m}(\CK_{i})$, in order to write 

\begin{equation}\label{tutto}
\begin{split}
&\ptm\left( -x<\sum_{j=1}^{\tnu}|\tbe_{j,\tnu}|(\tb_j U_j +(1-\tb_j)Z_j)<x\right)\\
&\qquad\leq \dfrac{1}{\a+1}\frac{k_m(a_m+\eta)}{A_m^{\a+1}K_{2,m}}+\ptm\{\max_{j=1,\dots,\tnu}|\tbe_{j,\tnu}(1-\tb_j)|>\veps\}\\
& \qquad+ \dfrac{4}{\pi}\int_0^{+\infty} \bar{H}_m(\xi,e^{-t\sigma})\xi^{\alpha-1}e^{-(a_m-\eta)\xi^{\alpha}}d\xi+\dfrac{16}{5\pi}M_m^2N_{2\alpha,m}e^{-t(1-2R_4)}\\
&\qquad+\dfrac{128}{25\pi}M_m^4N_{4\alpha,m} e^{-t(1-2R_6)}+\dfrac{2\textbf{c}\left \|g_{m}\right \|}{\tilde{d}_m}\left( e^{-\sigma t}+e^{-t(1-q\sigma\alpha /2-2R_q)}\right)\\
&\qquad+\dfrac{4}{\pi}\left \|v_m\right \|(N_{\alpha,m}+2N_{2\alpha,m}\left \|v_m\right \|)e^{-t(1-q\sigma\alpha /2-2R_q)}\\
&\qquad+\dfrac{2}{\pi}\left( x+A_m\right)\dfrac{1}{\alpha a_m^{1/\alpha}}\Gamma\left( \dfrac{1}{\alpha}\right). \end{split}
\end{equation}
\newline\newline

We set about the final step of the proof $-$ where the $i$-th summand in the RHS of \eqref{tutto} will be indicated by $\CR_i^{(m)}$ for $i=1,\dots,8$ $-$ whose aim is to associate suitable values of $m$ and $x$ with time $t$.\newline
\newline
\textbf{Step 3: }\textit{Quantification of $t\mapsto m(t)$ and $t\mapsto x(t)$.} As to $\CR_1^{m(t)}$, we show that \textit{any increasing function $t\mapsto m(t)$ which diverges, as $t\to+\infty$}, serves our purpose. Indeed, recalling the definitions of $k_m$ and $A_m$, 
\[
\CR_1^{(m(t))}=C_1 (1-F_0^*(x_{m(t)}))^\delta[1+r_1(t,m(t))]
\]
where $C_1:=2\dfrac{(1+2\pi)^{2(\alpha +1)}}{\alpha +1}$ and $r_1(t,m(t)):=\dfrac{\eta}{a_{m(t)}}$. 

We proceed to prove that \textit{$\CR_2^{(m(t))}$ vanishes at infinity whenever $m(\cdot)$ meets the condition
\begin{equation}\label{m'}
m(t)\leq \sup\left\lbrace n\in\mathbb{N}: x_n\leq (H_{2-\delta})^{-1}(e^{-t 2\tau_1/q})\right\rbrace\quad(t>0)
\end{equation}
$\tau_1$ being any point fixed in $(0,1-2R_q)$}. To this aim we note that, for $\veps=\veps_{m(t)}=A_{m(t)}\Big(\dfrac{3d^*}{8M_{m(t)}}\Big)^{1/\a}$,
\[
\CR_2^{(m(t))}\leq \ptm\{\max_{j=1,\dots,\tnu}|\tbe_{j,\tnu}|>\veps_{m(t)}\}\leq \veps_m^{-q/(1+p)}e^{-t(1-2R_q)}
\]
the last inequality being consequence of Lemma 1 in \cite{GabettaRegazzini2008}. Hence,

\[
\begin{split}
&\CR_2^{(m(t))}\\ 
&\; \leq \dfrac{\sqrt{8}^q(1+2\pi)^{q\a}}{(3d^*)^{q/2}}e^{-t(1-2R_q)}\Big(\dfrac{a_{m(t)}+\left \|v_{m(t)}\right \|}{k_{m(t)}^\a}\Big)^{q/2}\\
&\; \leq \dfrac{\sqrt{8}^q(1+2\pi)^{q\a}}{(3d^*)^{q/2}}e^{-t(1-2R_q)}\Big(\dfrac{2c_m+\a}{2^\a c_{m(t)} \Gamma(1-\a)\cos(\pi\a/2)}\dfrac{1}{(1-F^*_0(x_{m(t)}))^{2-\delta}}\\
&\;\times\Big[1+2^{\a}\dfrac{c_{m(t)}}{2c_{m(t)}+\a}(1-F^*_0(x_{m(t)}))\Big(\Gamma(1-\a)\cos(\pi\a/2)\\
&\qquad\qquad\qquad\qquad\qquad\qquad\qquad\qquad\qquad\qquad\qquad+\int_{(0,\pi)}(\sin x)/x^\a dx\Big)\Big]\Big)^{q/2}\\
& \;= C_2 \dfrac{e^{-t(1-2R_q)}}{(1-F^*_0(x_{m(t)}))^{q(2-\delta)/2}c_{m(t)}^{q/2}} [1+r_2(t,m(t))]
\end{split}
\]
with $C_2:= \dfrac{(1+2\pi)^{q\a}2^{(3-\a)q/2}\a^{q/2}}{(3d^*)^{q/2}\Gamma(1-\a)^{q/2}\cos^{q/2}(\pi\a/2)}$,  $\;r_2(t,m(t)):= \Big[1+2^{\a}(1-F^*_0(x_{m(t)}))$ $\Big(\Gamma(1-\a)\cos(\pi\a/2)+\int_{(0,\pi)}(\sin x)/x^\a dx\Big)\Big]^{q/2}-1$. Whence, in order that $\CR_2^{(m(t))}$ converge to $0$ as $t$ goes to infinity, it suffices that
\[
(H_{2-\delta}(x_{m(t)}))^{q/2}\geq e^{-t\tau_1}
\]
where $\tau_1$ is any point in the interval $(0,1-2R_q)$, a condition which is verified if
\[
x_{m(t)}\leq (H_{2-\delta})^{-1}(e^{-t 2\tau_1/q)})
\] 
and, \textit{a fortiori}, if $m(\cdot)$ meets \eqref{m'}. 

As far as $\CR_3^{(m(t))}$ is concerned, \textit{a suitable choice of $t\mapsto m(t)$, which is consistent with \eqref{m'}, is given by
\begin{equation}\label{m1}
m(t):= \sup\left\lbrace n\in\mathbb{N}:x_n\leq 2^{1/2}e^{\tau\sigma t/2}\right\rbrace
\end{equation}
where $\tau$ is a positive constant satisfying $2-\a-2\tau>0$.} This statement rests on Lemma \ref{g1m}, according to which
\begin{equation*}
\begin{split}
\CR_3^{(m(t))}&\leq \dfrac{4}{\pi}\Big(\int_0^{x_{m(t)}}+\int_{x_{m(t)}}^{+\infty}\Big)\sup_{0\leq s\leq e^{-\sigma t}}\left[ |\omega_{m(t)}(\xi s)|+4|\xi|^\alpha |\omega_{m(t)}(\xi s)|^2 \right]\\
&\qquad\qquad\qquad\qquad\qquad\times\xi^{\alpha-1}e^{-(a_{m(t)}-\eta)\xi^{\alpha}}d\xi\\
&+\dfrac{4}{\pi}\Big(\int_0^{x_{m(t)}}+\int_{x_{m(t)}}^{+\infty}\Big) \sup_{0\leq s\leq e^{-\sigma t}}\left[ |\rho_{m(t)}(\xi s)|+4|\xi|^\alpha |\rho_{m(t)}(\xi s)|^2 \right]\\
&\qquad\qquad\qquad\qquad\qquad\times\xi^{\alpha-1}e^{-(a_{m(t)}-\eta)\xi^{\alpha}}d\xi\\
&=: \Big(\CR_{3,1}^{(m(t))}+ \CR_{3,2}^{(m(t))}\Big)+\Big(\CR_{3,3}^{(m(t))}+ \CR_{3,4}^{(m(t))}\Big)
\end{split}
\end{equation*}
where $\omega_{m(t)}(\xi):=\dfrac{2}{\xi^\alpha}\int_0^{\xi x_{m(t)}} \sin x \Big[1-K_{1,m}\Big(F^*_0\Big(\dfrac{x}{\xi}\Big)+F^*_0(x_m)-1+\dfrac{\theta_m}{\a x_m^\a} \Big)\Big]dx$ and $\rho_{m(t)}(\xi)$ $:=\dfrac{2}{\a}K_{1,m}\theta_m\int_0^{\xi x_{m(t)}}\dfrac{\sin x}{x^\alpha}dx$ for every $\xi>0$. The analysis of $\CR_{3,1}^{(m(t))}$ is based on an argument, already used in the proof of Lemma \ref{g1m}, which yields $|\omega_{m(t)}(\xi)|\leq \min \left( \dfrac{x_{m(t)}^2 \xi^{2-\alpha}}{2},2\xi^{-\alpha}\right)$, and then
\begin{equation}\label{sup}
\sup_{0\leq s\leq e^{-\sigma t}}\left[ |\omega_{m(t)}(\xi s)|\right] \leq \left\lbrace 
\begin{array}{ll}
\displaystyle
\dfrac{x_{m(t)}^2 (\xi e^{-\sigma t})^{2-\alpha}}{2} & \mbox{ if } \xi e^{-\sigma t}\leq\dfrac{2}{x_{m(t)}}\\
2^{1-\alpha}x_{m(t)}^\alpha & \mbox{ if } \xi e^{-\sigma t}>\dfrac{2}{x_{m(t)}}.
\end{array}\right.
\end{equation}
The argument proceeds by assuming that $m=m(t)$ is the same as in \eqref{m1}. Hence, for every $\xi$ in $(0,x_{m(t)}]$, we get $x_{m(t)}\leq 2e^{\sigma t}/\xi$ and $\sup_{0\leq s\leq e^{-\sigma t}}$ $\left[ |\omega_{m(t)}(\xi s)|\right]\leq 2^{1-\alpha/2}\exp(-(2-\alpha-2\tau+\alpha \tau/2)\sigma t)$. In view of these inequalities,
\begin{equation*}\label{11}
\begin{split}
\CR_{3,1}^{(m(t))}& \leq\frac{2^{3-\alpha/2}}{\pi}\exp(-(2-\alpha-2\tau+\alpha \tau/2)\sigma t)\int_0^{x_{m(t)}}\xi^{\alpha-1}e^{-(a_{m(t)}-\eta)\xi^\alpha}d\xi\\
& +\frac{2^{6-\alpha}}{\pi}\exp(-2(2-\alpha-2\tau+\alpha \tau/2)\sigma t)\int_0^{x_{m(t)}}\xi^{2\alpha-1}e^{-(a_{m(t)}-\eta)\xi^\alpha}d\xi\\
& =\frac{2^{3-\alpha/2}}{\pi\a}\exp(-(2-\alpha-2\tau+\alpha \tau/2)\sigma t)\frac{1-e^{-(a_{m(t)}-\eta)x_{m(t)}^\alpha}}{(a_{m(t)}-\eta)}\\
& +\frac{2^{6-\alpha}}{\pi\a}\exp(-2(2-\alpha-2\tau+\alpha \tau/2)\sigma t)\left( \dfrac{1-e^{-(a_{m(t)}-\eta)x_{m(t)}^\alpha}}{(a_{m(t)}-\eta)^2}\right.\\
& \left. -\dfrac{x_{m(t)}^\alpha e^{-(a_{m(t)}-\eta)x_{m(t)}^\alpha}}{a_{m(t)}-\eta}\right) \\
&\leq C_{3,1} \dfrac{e^{-(2-\alpha-2\tau+\alpha \tau/2)\sigma t}}{a_{m(t)}-\eta}[1+r_{3,1}(t,m(t))]
\end{split}
\end{equation*}
where $C_{3,1}:= \dfrac{2^{3-\a/2}}{\pi\a}$ and $r_{3,1}(t,m(t)):=2^{3-\a/2} e^{-(2-\alpha-2\tau+\alpha \tau/2)\sigma t}/(a_{m(t)}-\eta)$. As to $\CR_{3,2}^{(m(t))}$, we resort to \eqref{sup} to write $\sup_{\xi>0,\;0\leq s\leq e^{-\sigma t}}\left[ |\omega_{m(t)}(\xi s)|\right]\leq 2^{1-\alpha}x_{m(t)}^\alpha$ and obtain
\begin{equation*}\label{12}
\begin{split}
\CR_{3,2}^{(m(t))}& \leq  C_{3,2} \dfrac{x_{m(t)}^{3\a}e^{-(a_{m(t)}-\eta)x_{m(t)}^\alpha}}{(a_{m(t)}-\eta)}[1+r_{3,2}(t,m(t))]
\end{split}
\end{equation*}
with $C_{3,2}:= \dfrac{2^{6-2\a}}{\pi\a}$ and $r_{3,2}(t,m(t)):=\dfrac{1}{2^{3-\a} x_{m(t)}^{2\a}}+\dfrac{1}{(a_{m(t)}-\eta) x_{m(t)}^{\a}} $. As to $\CR_{3,3}^{(m(t))}$, arguing in the same way as at the end of the proof of Lemma \ref{g1m}, we obtain
\begin{multline*}
\sup_{0\leq s\leq e^{-\sigma t}}\left[ \int_0^{x_{m(t)}\xi s}\dfrac{\sin x}{x^\alpha}dx \right] \\=\left\lbrace 
\begin{array}{ll}
\displaystyle
\int_0^{x_{m(t)}\xi e^{-\sigma t}}\dfrac{\sin x}{x^\alpha}dx\leq \dfrac{(x_{m(t)}\xi e^{-\sigma t})^{2-\alpha}}{2-\alpha} & \mbox{ if } x_{m(t)}\xi e^{-\sigma t}\leq \pi\\
\displaystyle
\int_{(0,\pi)}\dfrac{\sin x}{x^\alpha}dx\leq \dfrac{\pi^{2-\alpha}}{2-\alpha }& \mbox{ if } x_{m(t)}\xi e^{-\sigma t}> \pi.
\end{array}\right. 
\end{multline*}
Thanks to this inequality $-$ since \eqref{m1} with $\xi$ in $(0,x_{m(t)}]$ entails $x_{m(t)}\leq \pi e^{\sigma t}/\xi$ $-$ we have
\[
\sup_{0\leq s\leq e^{-\sigma t}}\left[ |\rho_{m(t)}(\xi s)| \right]\leq \dfrac{2^{3-\alpha/2}(1-F_0^*(x_{m(t)}))}{(2-\alpha)}\exp(-(2-\alpha-2\tau+\tau\alpha/2)\sigma t)
\]
which, in its turn, implies 
\begin{equation*}\label{21}
\CR_{3,3}^{(m(t))}\leq C_{3,3} \frac{e^{-(2-\alpha-2\tau+\tau\alpha/2)\sigma t}}{a_{m(t)}-\eta}[1+r_{3,3}(t,m(t))]
\end{equation*}
with $C_{3,3}:=\dfrac{2^{5-\a/2}}{\pi\a(2-\a)}$ and $r_{3,3}(t,m(t)):= \dfrac{2^{5-\a/2}}{2-\a} \frac{e^{-(2-\alpha-2\tau+\tau\alpha/2)\sigma t}}{a_{m(t)}-\eta}$. To deal with $\CR_{3,4}^{(m(t))}$, note that 
\[
\sup_{\xi>0,\;0\leq s\leq e^{-\sigma t}}\left[ |\rho_{m(t)}(\xi s)| \right]\leq \dfrac{(1-F_0^*(x_{m(t)}))x_{m(t)}^\alpha}{F_0^*(x_{m(t)})} \dfrac{\pi^{2-\alpha}}{2-\alpha }.
\]
Hence,
\begin{equation*}\label{22}
\CR_{3,4}^{(m(t))}\leq C_{3,4} \dfrac{x_{m(t)}^{3\a}e^{-(a_{m(t)}-\eta)x_{m(t)}^\alpha}}{a_{m(t)}-\eta}[1+r_{3,4}(t,m(t))]
\end{equation*}
with $C_{3,4}:= \dfrac{2^6\pi^{3-2\a}}{\a(2-\a)^2}$ and $r_{3,4}(t,m(t)):= \dfrac{2-\a}{2^3\pi^{2-\a}x_{m(t)}^{3\a}}+\dfrac{1}{x_{m(t)}^\a(a_{m(t)}-\eta)}$.

Proceeding, we now show that \textit{$\CR_4^{(m(t))}$ goes to $0$ as $t\to+\infty$, provided that
\begin{equation}\label{m2}
m(t)\leq \sup\left\lbrace n\in\mathbb{N}: x_n\leq (H_1)^{-1}(e^{-t\tau_2/2})\right\rbrace \qquad (t>0)
\end{equation}
with $\tau_2$ in $(0,1-2R_4)$}. It is worth noting that this condition is consistent with \eqref{m'}-\eqref{m1}. To verify this claim, recall that
\[
\begin{split}
N_{2\alpha,m(t)}&:=\int_{(0,+\infty)} e^{-(a_{m(t)}-\eta)\xi^\alpha}\xi^{2\alpha-1}d\xi=\dfrac{1}{\alpha(a_{m(t)}-\eta)^2}\\
&=\dfrac{1}{\alpha a_{m(t)}^2}\Big(\dfrac{a_{m(t)}}{a_{m(t)}-\eta}\Big)^2
\end{split}
\]
together with the definition of $M_{m(t)}$ and the bound to $\left\| v_m \right \|$ stated in Lemma \ref{g1m}. They yield
\[
\CR_4^{(m(t))}\leq C_4 \dfrac{e^{-t(1-2R_4)}}{((1-F^*_0(x_{m(t)})c_{m(t)})^2}[1+r_4(t,m(t))]
\]
with $C_4:= \dfrac{2^{4-2\a}3^2\a}{5\pi\Gamma^2(1-\a)\cos^2(\pi\a/2)}$ and $r_4(t,m(t)):= \Big(1+2^{\a}(1-F^*_0(x_{m(t)}))$ $\Big[\int_0^\pi \sin x/x^\a dx + \Gamma(1-\a)\cos(\pi\a/2)\Big]\Big)^2 \Big(\dfrac{a_{m(t)}}{a_{m(t)}-\eta}\Big)^2-1$. 

Analogously, we prove that \textit{$\CR_5^{(m(t))}$ vanishes at infinity if
\begin{equation}\label{m3}
m(t)\leq \sup\left\lbrace n\in\mathbb{N}: x_n\leq (H_1)^{-1}(e^{-t\tau_3/4})\right\rbrace
\end{equation}
with $\tau_3$ in $(0,1-2R_6)$}, a condition which is compatible with \eqref{m'},\eqref{m1} and \eqref{m2}. Indeed, integration by parts yields
\[
N_{4\alpha,m(t)}=\dfrac{6}{\alpha(a_{m(t)}-\eta)^4}=\dfrac{6}{\alpha a_{m(t)}^4}\Big(\dfrac{a_{m(t)}}{a_{m(t)}-\eta}\Big)^4
\]
and then, 
\begin{equation*}\label{R5}
\begin{split}
\CR_5^{(m(t))}\leq C_5 \dfrac{e^{-t(1-2R_6)}}{((1-F^*_0(x_{m(t)})c_{m(t)}))^4}[1+r_5(t,m(t))]
\end{split}
\end{equation*}
with $C_5:= \dfrac{3^5\cdot 2^{8-4\a}\a^3}{25\pi\Gamma^4(1-\a)\cos^4(\pi\a/2)}$ and $r_5(t,m(t)):=\Big(1+2^{\a}(1-F^*_0(x_{m(t)}))$ $\Big[\int_0^\pi \sin x/x^\a dx + \Gamma(1-\a)\cos(\pi\a/2)\Big]\Big)^4   \Big(\dfrac{a_{m(t)}}{a_{m(t)}-\eta}\Big)^4-1$. 

In order that $\CR_6^{(m(t))}$ vanishes at infinity, \textit{it is sufficient that  
\begin{equation}\label{m4}
m(t)\leq \sup\left\lbrace n\in\mathbb{N}: x_n\leq (H_1)^{-1}(e^{-t\tau_4\a})\right\rbrace
\end{equation}
with $\tau_4$ in $(0,\Lambda)$ and $\Lambda:=\min\lbrace \sigma,1-q\sigma\alpha /2-2R_q\rbrace$.} The proof of this statement is based on the classical inversion formula of c.f.'s, which allows us to write
\[
\begin{split}
\CR_6^{(m(t))}& = \dfrac{2\textbf{c}}{\tilde{d}_{m(t)}}\sup_{\xi\in\R}\Big[\dfrac{1}{2\pi}\int_\R e^{-iu\xi}e^{-a_{m(t)}|u|^\alpha}du\Big](e^{-\sigma t}+e^{-t(1-q\sigma\alpha /2-2R_q)})\\
&\leq \dfrac{2\textbf{c}8^{1/\a}\Gamma(1/\a)}{\pi(3d^*)^{1/\a}\a}e^{-\Lambda t}\Big(1+\dfrac{\left \|v_{m(t)}\right \|}{a_{m(t)}}\Big)^{1/\a}\\
& = C_6 \dfrac{e^{-\Lambda t}}{((1-F^*_0(x_{m(t)})c_{m(t)}))^{1/\a}}\Big[1+r_6(t,m(t))\Big]
\end{split}
\]
with $C_6:= \dfrac{2^{3/\a}\textbf{c}\Gamma(1/\a)\a^{1/\a}}{\Gamma^{1/\a}(1-\a)\cos^{1/\a}(\pi\a/2)\pi\a(d^*)^{1/\a}}$ and $r_6(t,m(t)):= \Big( 1+2^{\a}(1-F^*_0(x_{m(t)}))\Big[\Gamma(1-\a)\cos(\pi\a/2)+\int_0^\pi \sin x/x^\a dx\Big]\Big)^{1/\a}-1$. 

We next prove that \textit{$\CR_7^{(m(t))}$ vanishes at infinity if
\begin{equation}\label{m6}
m(t)\leq \sup\left\lbrace n\in\mathbb{N}: x_n\leq (H_1)^{-1}(1-e^{-t\tau_5/2})\right\rbrace
\end{equation} 
with $\tau_5$ in $(0,1-q\sigma\alpha /2-2R_q)$.} Indeed,
\[
\begin{split}
\CR_7^{(m(t))}&=\dfrac{4}{\pi}||v_{m(t)}||(N_{\alpha,m(t)}+2N_{2\alpha,m(t)}||v_{m(t)}||)e^{-t(1-q\sigma\alpha /2-2R_q)}\\
&=\dfrac{4}{\pi}||v_{m(t)}||N_{\alpha,m(t)}e^{-t(1-q\sigma\alpha /2-2R_q)}\\
&+\dfrac{8}{\pi}N_{2\alpha,m(t)}||v_{m(t)}||^2e^{-t(1-q\sigma\alpha /2-2R_q)}\\
&=:\CR_{7,1}^{(m(t))}+\CR_{7,2}^{(m(t))}
\end{split}
\]
and, by arguing as in similar previous cases,
\[
\CR_{7,1}^{(m(t))}\leq C_{7,1} \dfrac{e^{-t(1-q\sigma\alpha /2-2R_q)}}{(1-F_0^*(x_{m(t)}))c_{m(t)}}[1+r_{7,1}(t,m(t))]
\]
with $C_{7,1}:=\dfrac{3\cdot 2^{2-\a}}{\pi\Gamma(1-\a)\cos(\pi\a/2)}$, $r_{7,1}(t,m(t)):= \Big(1+ 2^{\a}(1-F^*_0(x_{m(t)}))$ $\times\int_{(0,\pi)}(\sin x)/x^\a dx\Big)\Big(\dfrac{a_{m(t)}}{a_{m(t)}-\eta}\Big)-1$, and 
\[
\CR_{7,2}^{(m(t))}\leq C_{7,2} \dfrac{e^{-t(1-q\sigma\alpha /2-2R_q)}}{((1-F_0^*(x_{m(t)}))c_{m(t)})^2}[1+r_{7,2}(t,m(t))]
\]
with $C_{7,2}:=\dfrac{3^2 2^{3-2\a}\a}{\pi\Gamma^2(1-\a)\cos^2(\pi\a/2)}$, $r_{7,2}(t,m(t)):= \Big(1+ 2^{\a}(1-F^*_0(x_{m(t)}))$ $\times\int_{(0,\pi)}(\sin x)/x^\a dx\Big)^2\Big(\dfrac{a_{m(t)}}{a_{m(t)}-\eta}\Big)^2-1$. 

Finally, 
\[
\CR_8^{(m(t))}=\dfrac{2x}{\pi\alpha a_{m(t)}^{1/\alpha}}\Gamma\left( \dfrac{1}{\alpha}\right)+\dfrac{2}{\pi}\dfrac{(1-F_0^*(x_{m(t)}))^{\frac{1-\delta}{\alpha}}}{\alpha(1+2\pi)^2}\Gamma\left( \dfrac{1}{\alpha}\right)
\]
which goes to $0$ whenever $m(t)$ diverges, as $t\to+\infty$.

In view of \eqref{m'}, \eqref{m1}, \eqref{m2}, \eqref{m3}, \eqref{m4}, \eqref{m6}, in conjunction with the fact that both $\CR_1^{(m(t))}$ and $\CR_8^{(m(t))}$ do not require specific bounds to $m(\cdot)$, an admissible form for it is
\[
m(t):=\sup \left\lbrace n\in\mathbb{N} : x_n\leq 2^{1/2}e^{\tau' t}\wedge(H_{1})^{-1}(e^{-t\tau_1''})\wedge(H_{2-\delta})^{-1}(e^{-t\tau_2''})\right\rbrace
\]
where $\tau':=\tau\sigma/2$, $\tau_1'':=\min\left\lbrace\dfrac{\tau_2}{2},\dfrac{\tau_3}{4},\a\tau_4,\dfrac{\tau_5}{2}\right\rbrace$, $\tau''_2:=2\tau_1/q$,  $\sigma$ and $q$ are fixed positive numbers and $\delta$ is a fixed element of $(0,1)$. 

To complete the argument, it must be observed that $x=x(t)$ appears only in $\CR_8^{(m(t))}$ which, to vanish at infinity, requires that $x=x(t)=a_{m(t)}^{1/\a}\veps(t)$, $\veps(\cdot)$ being a positive function vanishing at infinity. The function $x(\cdot)$ provides the desired lower bound to the rate of explosion.

\appendix
\section{Proof of Proposition \ref{beta}}
By virtue of the Skorokhod construction there is a subset $\cO'$ of $\cO$, with $\cP$-probability $1$, such that the recursive relations \eqref{ricorrenzabeta} hold true for every point of $\cO'$. Then, we redefine the $\cben$'s according to \eqref{ricorrenzabeta} outside $\cO'$. This does not change the distribution of $\cW_n$. Next, we introduce the compact space
\[
M:= \overline{\N}^\infty\times\Big(\times_{j\geq1}\overline{\N}_j^\infty\Big)\times\Big(\times_{j\geq1}I_j^\infty\Big)
\]
where $\overline{\N}_1,\overline{\N}_2,\dots$ are copies of $\N:=\{1,2,\dots,+\infty\}$ and $I_1,I_2,\dots$ are copies of $[0,2\pi]$. We define the vector-valued function $\widehat{Y}\colon \cO\to M$ such that 
\[
\widehat{Y}:=\Big( (\cnu^{(j)})_{j\geq1}, (\ci^{(j)})_{j\geq1}, (\cte^{(j)})_{j\geq1}\Big) 
\]
(same notation as in Section \ref{sec:preliminaries}), and the mappings
\[
\begin{split}
f_i(\widehat{Y})&:=\Big( (\cnu^{(1)},\ci^{(1)}),\dots,(\cnu^{(i)},\ci^{(i)}),\big(c_p(\cte_j^{(1)}),s_p(\cte_j^{(1)})\big)_{j\geq1},\qquad\qquad(\widehat{Y}\in M)\\
&\qquad\qquad\qquad\dots,\big(c_p(\cte_j^{(i)}),s_p(\cte_j^{(i)})\big)_{j\geq1}\Big)
\end{split}
\] 
$i=1,2,\dots$. It should be noted that $(\cnu^{(i)},\ci^{(i)})$ specifies a McKean tree, say $\cCT_i$. The leading idea of the proof is the construction of a decreasing sequence $(A_n)_{n\geq1}$ of nonempty closed subsets of $M$ such that inequalities \eqref{limB}-\eqref{sumB} are met simultaneously when $\widehat{Y}$ belongs to $A_n$, for every $n$. To show this we first exhibit, for each $n$, a distinguished McKean tree $\CT_n$ together with a specific sequence of angles $\theta^{(n)}$ for which the desired inequalities occur. Then, we will make use of the distributional properties of $(\cCT_n,\cten)$ to state the existence of suitable neighbourhoods of the above pair on which the inequalities of interest \eqref{limB}-\eqref{sumB} are preserved. The rule we follow to construct the $A_n$'s is recursive. In the sequel, we confine ourselves to describing the first two steps since the step from $A_n$ to $A_{n+1}$ is essentially the same but with a more complicated notation.\newline
\textit{Proof for $n=1$} Our aim is to combine a tree with a sequence of angles in such a way that they yield \eqref{limB}-\eqref{sumB}. We begin by focusing on the complete tree of suitable depth $N_1$ and on angles equal to $\pi/4$. The common value of the associated $\cbe$'s is $[(1/\sqrt{2})^{2/\a}]^{N_1}=(1/2)^{N_1/\a}$. We choose $N_1$ and $m_1$ in such a way that these $\cbe$'s are all greater than $1/(x_{m_1}-\veps)$, which is equivalent to $N_1\leq\log_2(x_{m_1}-\veps)^\a$. This inequality is satisfied if we select $m_1$ such that $x_{m_1}$ is greater than $(1+\veps)$ and $N_1=\lfloor \log_2(x_{m_1}-\veps)^\a\rfloor$. Now consider this complete McKean tree with $2^{N_1}$ leaves (see the construction of a generic tree at the end of Section 2) and assume all its leaves germinate to obtain the complete tree with $\nu^{(1)}=2^{N_1+1}$ leaves, that is $\CT_1$. According to the usual left-to-right order, $\cte^{(1)}_{2^{N_1}+k-1}$ will indicate the angle associated with the germination of the $k$-th leaf of the original tree (i.e. before germination), $k=1,\dots,2^{N_1}$. With a view to the construction of the $\cbe^{(1)}$'s, we start from a reference situation with $(2^{N_1}-1)$ angles $\theta^{(1)}_1=\dots=\theta^{(1)}_{2^{N_1}-1}=\pi/4$ and in which we choose values $\theta^{(1)}_{2^{N_1}+k-1}$, $k=1,\dots,2^{N_1}$, meeting $1/x_{m_1}\leq |\cos \theta^{(1)}_{2^{N_1}+k-1}|^{2/\a}(1/2)^{N_1/\a}\leq 1/(x_{m_1}-\veps)$. Values of this kind exist since $2^{N_1/\a}/(x_{m_1}-\veps)\leq 1$. With this choice of $(2^{N_1+1}-1)$ angles, inequality \eqref{limB} follows immediately. Moreover, it is easy to check that $\sum_{i=1}^{\lfloor \frac{\nu^{(1)}+1}{2}\rfloor}|\beta^{(1)}_{2i-1,\nu^{(1)}}|^\a\geq \frac{1}{x_{m_1}^\a}\lfloor\frac{\nu^{(1)}+1}{2}\rfloor\geq \frac{1}{2}(\frac{x_{m_1}-\veps}{x_{m_1}})^\a =:a>0$. Then, there are a tree, i.e. the above $\CT_1$ with $\nu^{(1)}=2^{N_1+1}$ leaves, and a nondegenerate closed interval $\CI_1$ including $(1/2\dots,1/2,|\cos \theta^{(1)}_{2^{N_1}}|^2,\dots,|\cos \theta^{(1)}_{2^{N_1+1}-1}|^2)$ with the following property: the $\cbe^{(1)}$'s associated with the values of $\Big((\cnu^{(1)},\ci^{(1)}),\; \big(c_p(\cte_j^{(1)}),s_p(\cte_j^{(1)})\big)_{j\geq1}\Big)$ contained in $\{\CT_1\}\times\CI_1\times[0,1]^\infty$ satisfy \eqref{limB}-\eqref{sumB}. It is of paramount importance to note that, thanks to the distributional properties of $(\cCT_1,\cte^{(1)})$, the probability of the event $\{ f_1(\widehat{Y})\in \{\CT_1\}\times\CI_1\times[0,1]^\infty\}$ is strictly positive. Now, $A_1:=f_1^{-1}\Big(\{\CT_1\}\times\CI_1\times[0,1]^\infty\Big)$ is a closed subset of $M$ thanks to the continuity of $f_1$.\newline
\textit{Proof for $n=2$.} We consider the above tree $\CT_1$, with $\nu^{(1)}$ leaves, in combination with the sequence $(\theta^{(1)}_1,\dots,\theta^{(1)}_{\nu^{(1)}})$ of angles associated with $\CT_1$ in the previous part of the proof. We next append a suitable complete tree to each leaf of $\CT_1$, in accord with the following rule, to obtain $\CT_2$. For the actual construction of the latter, we consider each node $k=1,\dots,2^{N_1}$ of depth $N_1$ in $\CT_1$ , and we replace its "left child" ("right child", respectively) with a complete tree of suitable depth $N^{(k)}_{2,1}$ ($N^{(k)}_{2,2}$, respectively). Arguing as in the first part of the proof, we determine $m_2$, $N^{(k)}_{2,1}$ and $N^{(k)}_{2,2}$ in such a way that $[(1/\sqrt{2})^{2/\a}]^{(N_1+N^{(k)}_{2,1})}|\cos \theta^{(1)}_{2^{N_1}+k-1}|^{2/\a}\geq 1/(x_{m_2}-\veps)$ and $[(1/\sqrt{2})^{2/\a}]^{(N_1+N^{(k)}_{2,2})}|\sin \theta^{(1)}_{2^{N_1}+k-1}|^{2/\a}\geq 1/(x_{m_2}-\veps)$, which are equivalent to $N^{(k)}_{2,1}\leq \log_2(x_{m_2}-\veps)^\a+\log_2|\cos \theta^{(1)}_{2^{N_1}+k-1}|^2-N_1$ and $N^{(k)}_{2,2}\leq \log_2(x_{m_2}-\veps)^\a+\log_2|\sin \theta^{(1)}_{2^{N_1}+k-1}|^2 -N_1$. Then, we define $m_2$ to be the first index $m>m_1$ an $x_m$ satisfies $\lfloor \log_2(x_m-\veps)^\a+\log_2|\beta^{(1)}_{j,\nu^{(1)}}|^\a \rfloor\geq1$ for every $j=1,\dots,\nu^{(1)}$. Then, we put $N^{(k)}_{2,1}:=\lfloor\log_2(x_{m_2}-\veps)^\a+\log_2|\cos \theta^{(1)}_{2^{N_1}+k-1}|^2 \rfloor -N_1$ and  $N^{(k)}_{2,2}:=\lfloor\log_2(x_{m_2}-\veps)^\a+\log_2|\sin \theta^{(1)}_{2^{N_1}+k-1}|^2 \rfloor -N_1$ for $k=1,\dots,2^{N_1}$. As to the remaining part of the argument, for the sake of notational simplicity we confine ourselves to giving a detailed description for $k=1$. We assume that each leaf of the complete tree of depth $N^{(1)}_{2,1}$, previously appended to the "left child" of the node $k=1$, germinates. We choose $\theta^{(2)}_{2^{N_1+1}+2^{N^{(1)}_{2,1}}+j-2}$ $-$ the angle associated with the $j$-th germination ($j=1,\dots,2^{N^{(1)}_{2,1}}-1$) $-$ in such a way that $1/x_{m_2}\leq (1/2)^{(N_1+N^{(1)}_{2,1})/\a}|\cos \theta^{(1)}_{2^{N_1}}|^{2/\a}|\cos \theta^{(2)}_{2^{N_1+1}+2^{N^{(1)}_{2,1}}+j-2}|^{2/\a}\leq 1/(x_{m_2}-\veps)$. As to its existence, it suffices to note that $2^{(N_1+N^{(1)}_{2,1})/\a}/[(x_{m_2}-\veps)|\cos \theta^{(1)}_{2^{N_1}}|^{2/\a}]\leq1$. We now repeat the procedure for the "right child" of the node $k=1$ with $N^{(1)}_{2,2}$ in place of $N^{(1)}_{2,1}$ and $\sin \theta^{(1)}_{2^{N_1}}$ in place of $\cos \theta^{(1)}_{2^{N_1}}$. The argument to extend  this construction to all the remaining nodes $k=2,\dots,2^{N_1}$ is essentially the same and requires only obvious small changes in notation. The resulting tree $-$ with $\nu^{(2)}=2\sum_{k=1}^{2^{N_1}}(2^{N^{(k)}_{2,1}}+2^{N^{(k)}_{2,2}})$ leaves $-$ is $\CT_2$. The angles to be associated with $\CT_2$ are $\theta^{(2)}_h=\theta^{(1)}_h$ for $h=1,\dots,2^{N_1+1}-1$, while every angle pertaining to the complete trees appended, of depths $N^{(\cdot)}_{2,1}$ and $N^{(\cdot)}_{2,2}$ respectively, is equal to $\pi/4$. It is now easy to show that $\beta^{(2)}$, generated by $\CT_2$ together with the aforesaid angles, satisfies \eqref{limB}-\eqref{sumB} with $K=2$. Moreover, the complete subtree $\CT'_2$ of $\CT_2$ having $\nu^{(1)}=2^{N_1}$ leaves coincides with $\CT_1$ and, as to the angles associated with this subtree, one has $\theta^{(2)}_h=\theta^{(1)}_h$  for $h=1,\dots,2^{N_1+1}-1$. As a consequence, \eqref{limB}-\eqref{sumB} are satisfied also for $K=1$. At this stage, we can state the existence of a tree, the same $\CT_2$ as above, and of a nondegenerate closed interval $\CI_2$ such that: (a) $\CI_2=\CI'_1\times\CI_{1,2}$ for which $\CI'_1\subset\CI_1$; (b) the $\cbe^{(2)}$'s associated with the values of $f_2(\widehat{Y})$ contained in $\{\CT_1\}\times\{\CT_2\}\times\CI_1\times[0,1]^\infty\times\CI_2\times[0,1]^\infty$ satisfy \eqref{limB}-\eqref{sumB}. Moreover, the event $\Big\{f_2(\widehat{Y})\in \{\CT_1\}\times\{\CT_2\}\times\CI_1\times[0,1]^\infty\times\CI_2\times[0,1]^\infty\Big\}$ has strictly positive probability. The set $A_2:=f_2^{-1}\Big(\{\CT_1\}\times\{\CT_2\}\times\CI_1\times[0,1]^\infty\times\CI_2\times[0,1]^\infty\Big)$ is a closed subset of $M$ and it is easy to conclude, from the definitions of $f_i$ ($i=1,2$), that $A_2\subset A_1$.\newline
\textit{Conclusion.} Assuming that the procedure has been repeated in the same way to obtain nonempty closed sets of $M\supset A_1\supset A_2\supset\dots\supset A_n\supset\dots$, the finite intersection principle leads us to conclude that $\bigcap_{n\geq1}A_n$ is nonempty. To complete the argument, it suffices to note that $\widehat{Y}^{-1}\Big(\bigcap_{n\geq1}A_n\Big)\neq\emptyset$ and \eqref{limB}-\eqref{sumB} are satisfied for every $\o^*$ in $\widehat{Y}^{-1}\Big(\bigcap_{n\geq1}A_n\Big)$.

\section{Proof of Proposition \ref{prop1}}
Considering the Skorokhod representation, in order that the law of \eqref{V_t},  under $\CP_t$, converges weakly, it is necessary that $\cLmn(\o)$ converges weakly as $n\to+\infty$, for every $\o$ in $\cO$. Then, assuming such a convergence, from the central limit theorem for symmetric triangular arrays (see Theorem 24 in Section 16.8 of \cite{FristedtGray}), for every $\o$ in $\cO$ there exists a L\'evy measure $\nu(\o)$ such that 
\begin{equation}\label{clt}
\nu(\o)[x,+\infty)=\lim_{n\to+\infty}\sum_{j=1}^{\cnun(\o)}\clmn_j(\o)[x,+\infty)
\end{equation}
for every positive $x$ for which $\nu(\o)\{x\}=0$. If \eqref{rho} is violated, there is $(x_m)_{m\geq1}$ satisfying \eqref{cond1} and, in view of Proposition \ref{beta}, there are $(m_n)_{n\geq1}$ and $\o^*$ in $\cO$ for which \eqref{limB}-\eqref{sumB} are valid for some $\veps>0$. There is also a strictly positive $x_0<1$ for which $\nu(\o^*)\{x_0\}=0$, and \eqref{clt}, with $\o=\o^*$ and $x=x_0$, becomes
\[
\begin{split}
&\nu(\o^*)[x_0,+\infty)\\
&\qquad=\lim_{n\to+\infty}\sum_{j=1}^{\cnun(\o^*)}\clmn_j(\o^*)[x_0,+\infty)\\
&\qquad= \lim_{n\to+\infty}\sum_{j=1}^{\cnun(\o^*)}\Big[1-F^*_0\Big(\frac{x_0}{|\cben_{j,\cnun(\o^*)}(\o^*)|}\Big)\Big]\\
&\qquad\geq \limsup_{n\to+\infty}\sum_{j=1}^{\cnun(\o^*)}\Big[1-F^*_0\Big(\frac{1}{|\cben_{j,\cnun(\o^*)}(\o^*)|}\Big)\Big]\quad\text{($x_0<1$)}\\
&\qquad \geq \limsup_{n\to+\infty}\sum_{i=1}^{\lfloor \frac{\cnun(\o^*)+1}{2}\rfloor}\Big[1-F^*_0\Big(\frac{1}{|\cben_{2i-1,\cnun(\o^*)}(\o^*)|}\Big)\Big]\frac{|\cben_{2i-1,\cnun(\o^*)}(\o^*)|^\a}{|\cben_{2i-1,\cnun(\o^*)}(\o^*)|^\a}\\
&\qquad \geq \limsup_{n\to+\infty}\sum_{i=1}^{\lfloor \frac{\cnun(\o^*)+1}{2}\rfloor}[1-F^*_0(x_{m_n})](x_{m_n}-\veps)^\a|\cben_{2i-1,\cnun(\o^*)}(\o^*)|^\a\\
&\qquad\qquad\qquad\qquad\qquad\qquad\qquad\qquad\qquad\qquad\qquad\qquad\qquad\text{(from \eqref{limB})}\\
&\qquad \geq a \limsup_{n\to+\infty} x_{m_n}^\a[1-F^*_0(x_{m_n})]\Big(\frac{x_{m_n}-\veps}{x_{m_n}}\Big)^\a \quad\text{(from \eqref{sumB})}\\
&\qquad= a \limsup_{n\to+\infty} \rho(x_{m_n}).
\end{split}
\]
This is a contradiction when \eqref{rho} is violated.

\section{Proof of Proposition \ref{prop2}}
In view of \eqref{V_t}, for every $B$ in $\CB(\R)$ we have
\[
\begin{split}
\mu(B,t)&=\CP_t(V\in B)\\
		&=\E_t\Big(\CP_t(V\in B\Big|\tnu)\Big)\\
		&= \sum_{n\geq1}\CP_t\{\tnu=n\}\CP_t\Big\{\sum_{j=1}^n\tbe_{j,n}X_j\in B\Big\}
\end{split}
\]
where the last equality holds since $\tnu$ is stochastically independent of $(X,\ti,\tte)$. Thus, to prove the proposition, it is enough to show that
\[
\sup_{n\geq1} \CP_t\Big\{\sum_{j=1}^n\tbe_{j,n}X_j\in I^c_\veps\Big\}\leq\veps
\]
holds for every $\veps>0$ and for a suitable interval $I_\veps$ such as that introduced immediately before the wording of the proposition. The proof is based on inequality \eqref{doob} applied to $Y=\sum_{j=1}^n\tbe_{j,n}X_j$, which has the same p.d. as $\sum_{j=1}^n|\tbe_{j,n}|X_j$, thanks to the symmetry of the common distribution $F^*_0$ of the $X$'s. The c.f. of the sum with absolute values of the coefficients, in view of the independence of $X$ and $\tbe$ and the mutual independence of the $X$'s, is given by $u\mapsto\E_t[\prod_{j=1}^n \vphi_0(u|\tbe_{j,n}|)]$, and hence
\begin{equation}\label{tight1}
\begin{split}
\CP_t\Big\{\Big|\sum_{j=1}^n|\tbe_{j,n}|X_j\Big|>C\Big\}&\leq\frac{1}{\Delta}\Big(1+\frac{2\pi}{C\Delta}\Big)^2\int_0^\Delta\Big[1-\E_t\Big(\prod_{j=1}^n \vphi^*_0(u|\tbe_{j,n}|)\Big)\Big]du\\
&= \E_t \Big(\frac{1}{\Delta}\Big(1+\frac{2\pi}{C\Delta}\Big)^2\int_0^\Delta\Big[1-\prod_{j=1}^n \vphi^*_0(u|\tbe_{j,n}|)\Big]du\Big).
\end{split}
\end{equation}
From now on, the proof proceeds, with slight changes, like in Step 4 in Section 3 of \cite{GabettaRegazzini2012}. If 
\begin{equation}\label{limiti}
0\leq 1-\vphi^*_0(u)\leq\frac{3}{8},
\end{equation} then there is $\theta$ such that $0\leq|\theta|\leq1$ and
\[
\log\vphi^*_0(u)=-\Big\{1-\vphi^*_0(u)\Big\}+\frac{4\theta}{5}\Big\{1-\vphi^*_0(u)\Big\}^2.
\]
See Lemma 3 of Section 3 in \cite{ChowTeicher}. Obviously, there is $\Delta_0>0$ so that $(\ref{limiti})$ holds true if $|u|\leq\Delta\leq\Delta_0$. For any $u$ of this kind, since $|\tbe_{j,n}|\leq1$ for every $j$ and $n$, it turns out that $u|\tbe_{j,n}|\leq \Delta$ and 
\[
\log\vphi^*_0(u|\tbe_{j,n}|)=-\Big\{1-\vphi^*_0(u|\tbe_{j,n}|)\Big\}+\frac{4\theta}{5}\Big\{1-\vphi^*_0(u|\tbe_{j,n}|)\Big\}^2
\]
is valid for every $j$ and $n$. Thus, the last integrand in \eqref{tight1} can be developed as follows
\[
\begin{split}
1-\prod_{j=1}^n \vphi^*_0(u|\tbe_{j,n}|)&=1-\exp\Big(\sum_{j=1}^n\log\vphi^*_0(u|\tbe_{j,n}|)\Big)\\
& =1-\exp\Big(\sum_{j=1}^n\Big[-(1-\vphi^*_0(u|\tbe_{j,n}|))+\frac{4}{5}\theta(1-\vphi^*_0(u|\tbe_{j,n}|))^2 \Big]\Big)\\
&\leq 1-\exp\Big(\sum_{j=1}^n\Big[-(1-\vphi^*_0(u|\tbe_{j,n}|))-\frac{3}{10}(1-\vphi^*_0(u|\tbe_{j,n}|))\Big]\Big)\\
& = 1-\exp\Big(-\frac{13}{10}\sum_{j=1}^n(1-\vphi^*_0(u|\tbe_{j,n}|))\Big)\leq \frac{13}{10}\sum_{j=1}^n(1-\vphi^*_0(u|\tbe_{j,n}|)).
\end{split}
\]
This entails
\[
\CP_t\Big\{\Big|\sum_{j=1}^n|\tbe_{j,n}|X_j\Big|>C\Big\}\leq \frac{13}{10}\frac{1}{\Delta}\Big(1+\frac{2\pi}{C\Delta}\Big)^2\E_t\Big(\sum_{j=1}^n\int_0^\Delta\Big[1-\vphi^*_0(u|\tbe_{j,n}|)\Big]du\Big).
\]
Now, proceeding as in Step 4 of \cite{GabettaRegazzini2012}, with $|\tbe_{j,n}|$ in place of $1/n^{1/\a}$,  
\[
\begin{split}
&\int_0^\Delta\Big[1-\vphi^*_0(u|\tbe_{j,n}|)\Big]du\\
&\qquad\qquad\qquad\leq 2\Delta^{\a+1}|\tbe_{j,n}|^\a M \int_0^{+\infty}\Big(\I_{\{(0,\veps)\}}(u)\frac{y^{1-\a}}{6}+\I_{\{(\veps,+\infty)\}}(u)\frac{1+u}{u^{\a+2}}\Big)du
\end{split}
\]
where $M:=\sup_{x\geq0}\rho(x)$ which is finite by hypothesis. After noticing that the integral in the right-hand-side is finite, and denoting its value by $A$, we have
\[
\CP_t\Big\{\Big|\sum_{j=1}^n|\tbe_{j,n}|X_j\Big|>C\Big\}\leq \frac{13}{5}MA\Delta^\a\Big(1+\frac{2\pi}{C\Delta}\Big)^2
\]
and hence, taking $C:=C_\veps=\Delta^{-1}\geq \Big(\frac{13}{5}\frac{MA(1+2\pi)^2}{\veps}\Big)^{1/\a}$, 
\[
\CP_t\{\Big|\sum_{j=1}^n|\tbe_{j,n}|X_j\Big|>C\}\leq \veps
\] 
for every $n$, i.e. $\sup_{n\geq1} \CP_t\Big\{\sum_{j=1}^n|\tbe_{j,n}|X_j\in I^c_\veps\Big\}\leq\veps$ with $I_\veps:=[-C_\veps,C_\veps]$.

\section{Rates of explosion in Examples \ref{es2}} In the probabilistic setting at the beginning of Section \ref{sec:preliminaries}, let us replace each $\CP_t$ by $\bar\CP_t$, with the proviso that, under $\bar\CP_t$, all the random elements maintain the previous p.d.'s, with the exception of the $X_j$'s whose common p.d.f. becomes $\bar F^*_{0,t}(u):= F^*_0(u \cdot \xi(t))$ for every $u>0$ and for $\xi(t):=x(t)/\veps_1(t)$ where $t\mapsto \veps_1(t)$ is a strictly positive function on $(0,+\infty)$ such that $\veps_1(t)\searrow 0$, as $t\to+\infty$, and $t\mapsto x(t)$ is the same as the one specified in each of the examples we are dealing with. Then, for every $x>0$,
\[
\mu\Big((-x\cdot \xi(t), x\cdot\xi(t)),t\Big)=\bar\CP_t\Big\{\sum_{j=1}^{\tnu}|\tbe_{j,\tnu}|X_j\in (-x,x)\Big\}.
\]
It is clear that for every $\veps_1(\cdot)$ such that $\mu\Big((-x\cdot \xi(t), x\cdot\xi(t)),t\Big) \to 1$, as $t\to+\infty$, $\xi(t)$ represents an upper bound to the rate of explosion $a_t$. To find such a kind of $\veps_1(\cdot)$'s, we can resort to the degenerate convergence criterion provided by the central limit theorem (see, e.g. \cite{Loeve}). Accordingly, if
\begin{equation}\label{eq1tcl}
\bar\E_t\Big(\sum_{j=1}^{\tnu}\int_{|x|\geq\veps}d\bar F^*_{0,t}\Big(\dfrac{x}{|\tbe_{j,\tnu}|}\Big)\Big)\to 0
\end{equation}
and
\begin{equation}\label{eq2tcl}
\bar\E_t\Big(\sum_{j=1}^{\tnu}\int_{-\tau}^\tau x^2 d\bar F^*_{0,t}\Big(\dfrac{x}{|\tbe_{j,\tnu}|}\Big)\Big)\to 0
\end{equation}
hold, as $t\to+\infty$, for every $\veps>0$ and for some $\tau>0$, then there are a strictly increasing and divergent sequence of times $t_n$ and a version  of the conditional distribution of $\sum_{j=1}^{\tnu}|\tbe_{j,\tnu}|X_j$, given $(\tnu,\ti,\tte)$, which converges weakly to the unit mass $\delta_0$ at $0$. (It is worth noting that, for each $t_n$, the above-mentioned conditional distribution is derived assuming that the reference p.d. on $(\O,\CF)$ is just $\bar\CP_{t_n}$.) Then, $\mu\Big((-x\cdot \xi(t_n), x\cdot\xi(t_n)),t_n\Big)\to 1$ for every $x>0$. Since the same argument can be used to prove that every strictly increasing and divergent sequence of times $\tau_n$ includes a suitable subsequence $\tau'_n$ of the same type such that $\mu\Big((-x\cdot \xi(\tau'_n), x\cdot\xi(\tau'_n)),t_n\Big)\to1$, we get: If \eqref{eq1tcl}-\eqref{eq2tcl} hold, then $\mu\Big((-x\cdot \xi(t), x\cdot\xi(t)),t\Big)\to 1$ for every $x>0$. 

Applying this criterion to the first case considered in Examples \ref{es2}, we have, for sufficiently large values of $t$,
\[
\bar\E_t\Big(\sum_{j=1}^{\tnu}\int_{|x|\geq\veps}d\bar F^*_{0,t}\Big(\dfrac{x}{|\tbe_{j,\tnu}|}\Big)\Big)=\frac{2}{\veps^\beta\xi(t)^\beta}\bar\E_t\sum_{j=1}^{\tnu}|\tbe_{j,\tnu}|^\beta=\frac{2}{\veps^\beta\xi(t)^\beta}\exp[t(2R_{2\beta/\a}-1)]
\]
where the last equality follows from Proposition 8 in \cite{GabReg2006}. Moreover, analogously,
\[
\bar\E_t\Big(\sum_{j=1}^{\tnu}\int_{-\tau}^\tau x^2 d\bar F^*_{0,t}\Big(\dfrac{x}{|\tbe_{j,\tnu}|}\Big)\Big)=\dfrac{\tau^{2-\beta}\beta}{(2-\beta)\xi(t)^\beta}\exp[t(2R_{2\beta/\a}-1)].
\]
At this stage, it is easy to verify that, for $\beta\to 0^+$ of $\beta\to\a^-$, in order that \eqref{eq1tcl}-\eqref{eq2tcl} be satisfied, it is necessary to take $\xi(t)=\dfrac{C'}{\veps_1(t)^{1/\beta}}e^{t[2R_{2\beta/\a}-1+B(\a-\beta)\beta\a]/\beta}$. Moreover, with such a choice of $\xi(\cdot)$, \eqref{eq1tcl}-\eqref{eq2tcl} turn out to be valid for every $\beta$ in $(0,\a)$. 

As far as the second case in Examples \ref{es2} is concerned, according to the same way of reasoning, we observe that
\[
\bar\E_t\Big(\sum_{j=1}^{\tnu}\int_{|x|\geq\veps}d\bar F^*_{0,t}\Big(\dfrac{x}{|\tbe_{j,\tnu}|}\Big)\Big)\leq \dfrac{2\log (\veps\xi(t))}{\veps^\a\xi(t)^\a}+\dfrac{2}{\veps^\a\xi(t)^\a}\bar\E_t\Big(\sum_{j=1}^{\tnu}|\tbe_{j,\tnu}|^\a\log\dfrac{1}{|\tbe_{j,\nu}|}\Big).
\]
The former summand, in the RHS of the above inequality, goes to $0$ as $t\to+\infty$. As to the latter summand,
\[
\begin{split}
&\dfrac{2}{\veps^\a\xi(t)^\a}\bar\E_t\Big(\sum_{j=1}^{\tnu}|\tbe_{j,\tnu}|^\a\log\dfrac{1}{|\tbe_{j,\nu}|}\Big)\\
&\quad \leq \dfrac{4}{\a\veps^\a\xi(t)^\a}\bar\E_t\Big(\log\sum_{j=1}^{\tnu}|\tbe_{j,\tnu}|^\a\dfrac{1}{|\tbe_{j,\nu}|^{\a/2}}\Big)\qquad\text{(in view of the Jensen inequality,}\\
&\qquad\qquad\qquad\qquad\qquad\qquad\qquad\qquad\qquad\qquad\quad \text{since $\sum_{j=1}^{\tnu}|\tbe_{j,\tnu}|^\a=1$)}\\
&\quad \leq \dfrac{4}{\a\veps^\a\xi(t)^\a}\log\bar\E_t\Big(\sum_{j=1}^{\tnu}|\tbe_{j,\tnu}|^{\a/2}\Big)\qquad\quad\qquad\text{(from the Jensen inequality)}\\
&\quad \leq \dfrac{4\veps_1(t)^\a}{\a\veps^\a c^\a}\Big(\dfrac{4}{\pi}-1\Big)\qquad\qquad\qquad\qquad\qquad\quad\text{(from Proposition 8 in \cite{GabReg2006})}
\end{split}
\] 
which completes the argument to prove \eqref{eq1tcl}. Finally, the validity of \eqref{eq2tcl} can be verified, after an integration by parts, in the same way.

\section{Proof of Lemma \ref{g1m}} The first statement is an immediate consequence of well-known properties of the domain of attraction described, for example, in Section 2.6 of \cite{IbragimovLinnik}. Accordingly, it turns out that $a_m=(\pi(c_1+c_2))/(2\Gamma(\a)\sin(\pi\a/2))$ where $c_1=c_2=\dfrac{K_{1,m}\theta_m}{\a}$ and hence
\[
\begin{split}
a_m=&\dfrac{K_{1,m}\theta_m}{\a}\frac{\pi}{\Gamma(\a)\sin(\pi\a/2)}\\
&=\dfrac{K_{1,m}\theta_m}{\a}\frac{\Gamma(1-\a)\sin(\pi\a)}{\sin(\pi\a/2)}\\
&\qquad\qquad\qquad\qquad\qquad\text{(by the reflection property of the $\Gamma$ function)}\\
&= \dfrac{2}{\a}K_{1,m}\theta_m \Gamma(1-\a)\cos(\pi\a/2)\\
&= \dfrac{2}{\a}K_{1,m}\theta_m\int_{(0,+\infty)}\frac{\sin x}{x^\a}dx \qquad\text{(by Euler's formula).}
\end{split}
\]
To prove that $1-\hat{g}_{1,m}(\xi)=(a_m+v_m(\xi))\xi^\alpha$, set $F_m(x):=1-G_{1,m}(x)=G_{1,m}(-x)$ ($x>0$) to write
\[
1-\hat{g}_{1,m}(\xi)=\int_0^{+\infty}(e^{i\xi x}-1)dF_m(x)+\int_0^{+\infty}(e^{-i\xi x}-1)dF_m(x).
\]
As to the first integral, 
\begin{equation}\label{n1}
\begin{split}
&\int_0^{+\infty}(e^{i\xi x}-1)dF_m(x)\\
&\quad= -i\xi\int_0^{+\infty}\Big(\cos(\xi x)+ i \sin(\xi x)\Big)F_m(x)dx\quad \text{(by integration by parts)}\\
&\quad= -i\xi^\alpha \int_{(0,+\infty)} \cos(x)\dfrac{h_1(x/\xi)}{x^\alpha}dx+\xi^\alpha\int_{(0,+\infty)} \sin(x)\dfrac{h_1(x/\xi)}{x^\alpha}dx
\end{split}
\end{equation}
where 
\begin{equation}\label{h1}
h_1(x):= 1-K_{1,m}\Big(F^*_0(x)+F^*_0(x_m)-1+\dfrac{\theta_m}{\a x_m^\a}\Big) x^\alpha\I_{\{0<x< x_m\}} + \dfrac{K_{1,m}\theta_m}{\a}\I_{\{x\geq x_m\}}.
\end{equation}
Analogously,
\begin{equation}\label{n2}
\begin{split}
&\int_0^{+\infty}(e^{-i\xi x}-1)dF_m(x)\\
&\qquad= i\xi^\alpha \int_{[0,+\infty)} \cos(x)\dfrac{h_1(x/\xi)}{x^\alpha}dx+\xi^\alpha\int_{[0,+\infty)} \sin(x)\dfrac{h_1(x/\xi)}{x^\alpha}dx.
\end{split}
\end{equation}
Combination of \eqref{n1}, \eqref{n2} and \eqref{h1} gives
\[
\dfrac{1-\hat{g}_{1,m}(\xi)}{\xi^\alpha}=a_m+v_m(\xi)
\]
where $a_m$ is the same as in \eqref{am} and $v_m$ is the same as in the wording of the lemma, i.e.
\[
v_m(\xi)=\dfrac{2}{\xi^\alpha}\int_0^{\xi x_m} \sin x\; \o_{m,\xi}(x)dx-\dfrac{2}{\a}K_{1,m}\theta_m\int_0^{\xi x_m}\dfrac{\sin x}{x^\alpha}dx
\]
where $\o_{m,\xi}:=1-K_{1,m}\Big[F^*_0\Big(\frac{x}{\xi}\Big)+F^*_0(x_m)-1+\dfrac{\theta_m}{\a x_m^\a}\Big)$ for $x<\xi x_m$. 
Since $0<\a<2$, we have $\int_0^{\xi x_m}(\sin x)/x^\a dx\to 0$ as $\xi\to 0^+$ and, then, the latter summand is $o(1)$ as $\xi\to 0^+$. As to the former, 
\[
\left| \dfrac{2}{\xi^\alpha}\int_0^{\xi x_m} \sin x \;\o_{m,\xi}(x)dx\right|\leq \dfrac{2}{\xi^\alpha}\int_0^{\xi x_m} |\sin x| \sup_{t\in[0,x_m]} \o_{m,\xi}(x)dx.
\]
Observe that $x\mapsto \o_{m,\xi}(x)$ is nonnegative, monotonically nonincreasing on $[0,\xi x_m]$ and its value at $0$ is less or equal to $1/2$. Then, for $\xi\leq \pi/x_m$
\[
\left| \dfrac{2}{\xi^\alpha}\int_0^{\xi x_m} \sin x \;\o_{m,\xi}(x) dx\right|\leq \dfrac{1}{\xi^\a}\int_0^{\xi x_m}\sin x dx= \dfrac{1}{2}x_m^2\xi^{2-\a}+ o(|\xi|^{3-\a})
\]
which completes the proof that $v_m$ is $o(1)$ as $\xi\to 0^+$. As to the supremum norm of $v_m$,
\begin{equation}\label{svm}
\begin{array}{l}
\displaystyle
\left\|v_m\right\|:=\sup_{\xi\in\R^+} |v_m(\xi)|\leq \sup_{\xi\in\R^+}\left|\dfrac{1}{\xi^\alpha}\int_0^{\xi x_m} \sin x \;\o_{m,\xi}(x)dx\right| +\\
\displaystyle
\;\;\;\;\;\;\;\;\;\;\;\;\;\;\;\;\;\;\;\;+ \sup_{\xi\in\R^+}\left|\dfrac{2}{\a}K_{1,m}\theta_m\int_0^{\xi x_m}\dfrac{\sin x}{x^\alpha}dx\right|. 
\end{array}
\end{equation}
By the second mean value theorem, there is $a$ in $[0,\xi x_m)$ such that
\[
\left|\int_0^{\xi x_m} \omega_{m,\xi}(x)\sin x dx\right|=\left| \omega_{m,\xi}(0)\int_0^a \sin xdx\right| \leq 2.
\]
Moreover, from $|\sin x|\leq x$ for any $x\geq0$,
\begin{equation*}
\left| \int_{0}^{\xi x_m} \omega_{m,\xi}(x)\sin xdx\right|\leq \omega_{m,\xi}(0)\int_{0}^{\xi x_m}|\sin x|dx\leq \dfrac{x_m^2\xi^2}{2}
\end{equation*}
and combination of the above bounds gives 
\[
\left|\dfrac{1}{\xi^\alpha} \int_{0}^{\xi x_m} \omega_{m,\xi}(x)\sin x\;dx\right|\leq \min\left\lbrace  \dfrac{x_m^2\xi^{2-\alpha}}{2},2\xi^{-\alpha}\right\rbrace\leq 2^{1-\a}x_m^\a.
\]
To bound the latter summand in the RHS of \eqref{svm}, we begin by recalling the obvious inequalities
\begin{equation*}
\dfrac{2}{((k+1)\pi)^\alpha}\leq \left| \int_{k\pi}^{(k+1)\pi} \dfrac{\sin x}{x^\alpha}dx\right|\leq \dfrac{2}{(k\pi)^\alpha}\qquad(k=1,2,\dots)
\end{equation*}
and by setting $k^*=k^*(\xi):=\max\{k\in\N:\; k\pi\leq\xi x_m\}$ for every $\xi>0$. These inequalities, combined with
\begin{multline*}
\int_0^{\xi x_m}\dfrac{\sin x}{x^\alpha}dx=\left|\int_0^{\pi}\dfrac{\sin x}{x^\alpha}dx\right|-\left|\int_{\pi}^{2\pi}\dfrac{\sin x}{x^\alpha}dx\right|+....\\
...+(-1)^{k^*-1}\left|\int_{(k^*-1)\pi}^{k^*\pi}\dfrac{\sin x}{x^\alpha}dx\right|+(-1)^{k^*}\left|\int_{k^*\pi}^{x_m\xi}\dfrac{\sin x}{x^\alpha}dx\right|
\end{multline*}
yield\newline
\begin{equation*}
\begin{split}
\int_0^{\xi x_m}\dfrac{\sin x}{x^\alpha}dx&\leq \int_0^\pi \dfrac{\sin x}{x^\alpha}dx-\dfrac{2}{(2\pi)^\alpha} +\dfrac{2}{(2\pi)^\alpha}+...-\dfrac{2}{(x_m\xi)^\alpha}\I_{\lbrace \text{odd numbers}\rbrace}(k^*)\\
&\leq \int_0^\pi \dfrac{\sin x}{x^\alpha}dx
\end{split}
\end{equation*}
and
\begin{equation*}
\int_0^{\xi x_m}\dfrac{\sin x}{x^\alpha}dx\geq \dfrac{2}{\pi^\alpha} -\dfrac{2}{\pi^\alpha}+...+\dfrac{2}{(x_m\xi)^\alpha}\I_{\lbrace \text{even numbers}\rbrace}(k^*)\geq 0.
\end{equation*}
\newline
This entails 
\[
\sup_{\xi\in\R^+}\left|\dfrac{2}{\a}K_{1,m}\theta_m\int_0^{\xi x_m}\dfrac{\sin x}{x^\alpha}dx\right|\leq \dfrac{2}{\a}K_{1,m}\theta_m\int_0^\pi \dfrac{\sin x}{x^\alpha}dx.
\]


\begin{acknowledgements}
We would like to thank two anonymous referees for giving several valuable comments regarding the presentation of the main result.
\end{acknowledgements}



\end{document}